\documentclass[letterpaper, 10 pt, conference]{ieeeconf} 

\IEEEoverridecommandlockouts                              

\usepackage{csquotes}
\usepackage[usenames, dvipsnames]{color}
\usepackage{algorithm2e}
\usepackage{amsmath, amssymb, amsfonts, subfig}
\usepackage{soul}
\usepackage{tikz}
\usetikzlibrary{positioning}
\usepackage{flushend}
\usepackage{environ}
\usepackage{array,booktabs,amsmath,amssymb,amsfonts,textgreek,graphicx,fancyvrb,tabularx, xcolor, soul}
\usepackage{mathtools}
\SetAlCapSkip{1em}
\newif\ifarticle \articletrue

\title{\LARGE \bf
Cybersecurity in Distributed and Fully-Decentralized Optimization: Distortions, Noise Injection, and ADMM*
}

\author{Eric Munsing$^{1}$ and Scott Moura$^{1,2}$
\thanks{*The work was
  supported by the National Science Foundation Graduate Research Fellowship under Grant No. 1106400}
\thanks{$^{1}$Eric Munsing and Scott Moura are with the Department of Civil and Environmental Engineering, University of California at Berkeley, USA {\tt\small e.munsing@berkeley.edu}}%
\thanks{$^{2}$Scott Moura is with the Tsinghua-Berkeley Shenzhen Institute, University of California at Berkeley, USA}
}
\begin{document}

\maketitle
\thispagestyle{empty}
\pagestyle{empty}

\begin{abstract}
As problems in machine learning, smartgrid dispatch, and IoT coordination problems have grown, \textit{distributed} and \textit{fully-decentralized} optimization models have gained attention for providing computational scalability to optimization tools.  However, in applications where consumer devices are trusted to serve as distributed computing nodes, compromised devices can expose the optimization algorithm to cybersecurity threats which have not been examined in previous literature. This paper examines potential attack vectors for generalized distributed optimization problems, with a focus on the \textit{Alternating Direction Method of Multipliers} (ADMM), a popular tool for convex optimization. Methods for detecting and mitigating attacks in ADMM problems are described, and simulations demonstrate the efficacy of the proposed models. The weaknesses of fully-decentralized optimization schemes, in which nodes communicate directly with neighbors, is demonstrated, and a number of potential architectures for providing security to these networks is discussed.
\end{abstract}

\section{Motivation}

With increased computing power, convex optimization tools have gathered attention as a fast method for solving constrained optimization problems. However, some problems are still too large to be solved centrally- e.g. scheduling the energy consumption of millions of electric vehicles. 

In these applications, \textit{decentralized optimization} techniques break the problem into a set of subproblems which can be rapidly solved on distributed computing resources, with an \textit{aggregator} or \textit{fusion node} bringing the distributed problems into consensus on a global solution. This distribution can yield significant computational benefits, greatly speeding computation times.

\begin{figure}
	\centering
	\ifarticle \newcommand{\figwidth}{1.0} \else \newcommand{\figwidth}{1.0} \fi
	\includegraphics[width=\figwidth\columnwidth]{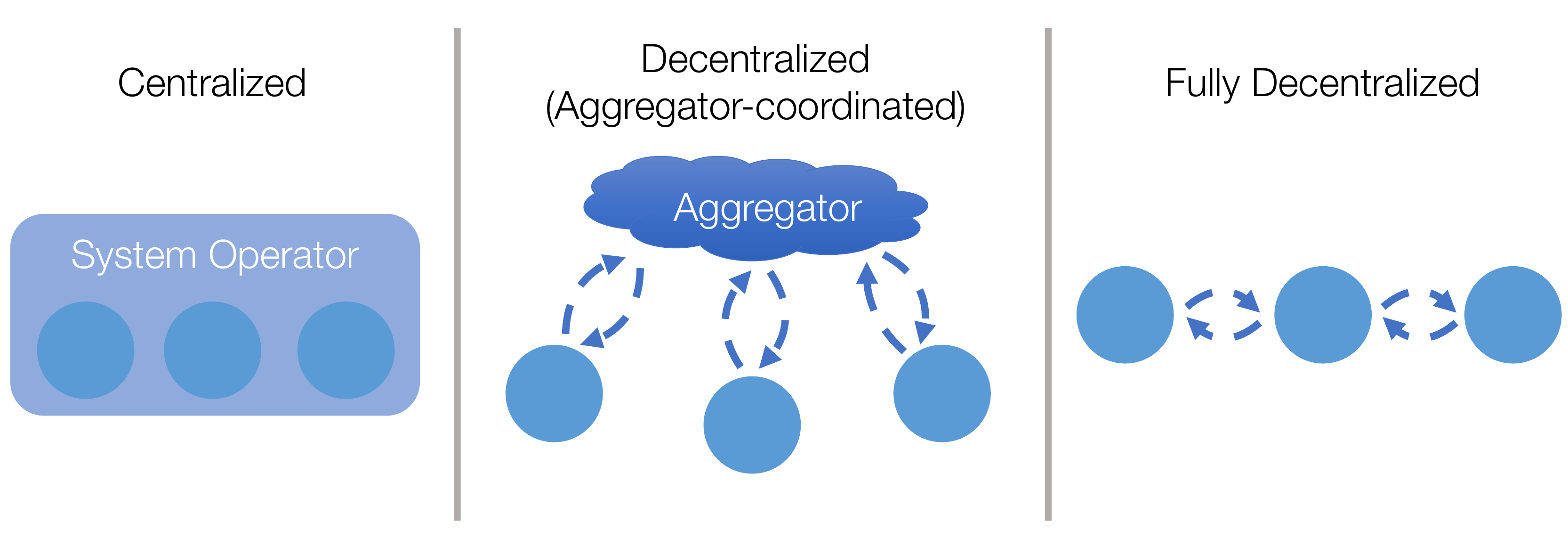}
	\caption{Different structures for solving mathematical optimization problems, from left: Centralized optimization, where all details of objective and constraints are held by a central entity. Decentralized optimization (also called aggregator-coordinated optimization) where local nodes hold local objective and constraint information, and an aggregator brings nodes into consensus on shared constraints. Fully-decentralized optimization, where no centralized entity exists but neighbors communicate directly with each other to achieve consensus.}
	\label{fig:optimizationtypes}
\end{figure}

A further development on this model, \textit{fully-decentralized} optimization techniques remove the aggregator, and instead let nodes reach consensus with neighbors, in a process which gradually brings the entire system into consensus (see Figure \ref{fig:optimizationtypes}). This approach significantly reduces communication requirements relative to aggregator-coordinated systems.

The potential for reducing computation time and communication burden makes these distributed and fully-decentralized models particularly compelling in smartgrid applications, where systems may need to scale to millions of devices.  However, by moving the computation from enterprise servers to consumer devices, significant new security weaknesses are exposed.  This article outlines architectures and algorithms for addressing these weaknesses.

\section{Background Literature}

The vulnerability of physical infrastructure to cyberattacks was dramatically highlighted in 2000 when the SCADA system controlling the Maroochy Water Plant in Australia was remotely hijacked, resulting in the release of over one million gallons of untreated sewage \cite{10.1007/978-0-387-75462-8_6}. A decade later, the US government developed the \textit{Stuxnet} worm to attack programmable logic controllers which operated uranium enrichment centrifuges managed by the Iranian government \cite{Langner2011}; the Iranian government recently launched similar attacsk against targets in Saudi Arabia \cite{Greenberg2017a}. A descendant of the Stuxnet worm was used more dramatically by the Russian government to repeatedly cause widespread blackouts in the Ukranian electricity during 2016 and 2017 \cite{Greenberg2017, Conway2016}.  Recent incursions into the networks of American electricity operators by Russian hackers may be laying the foundation for future incursions against U.S. infrastructure \cite{Perlroth2018, Atherton2018}.

 Against this backdrop of cyberattacks used as national security threats, President Obama issued an Executive Order \cite{EO13636} and Presidential Policy Directive \cite{PPD21, ISC2015} directing the executive branch to improve cybersecurity of critical infrastructure. In the face of proven vulnerabilities in consumer-level smart grid devices \cite{Rouf2012} it seems critical that distributed algorithms for smartgrid controls be designed to be able to identify, localize, and mitigate the effects of these attacks.

\begin{figure}
	\centering
	\includegraphics[scale = 0.5]{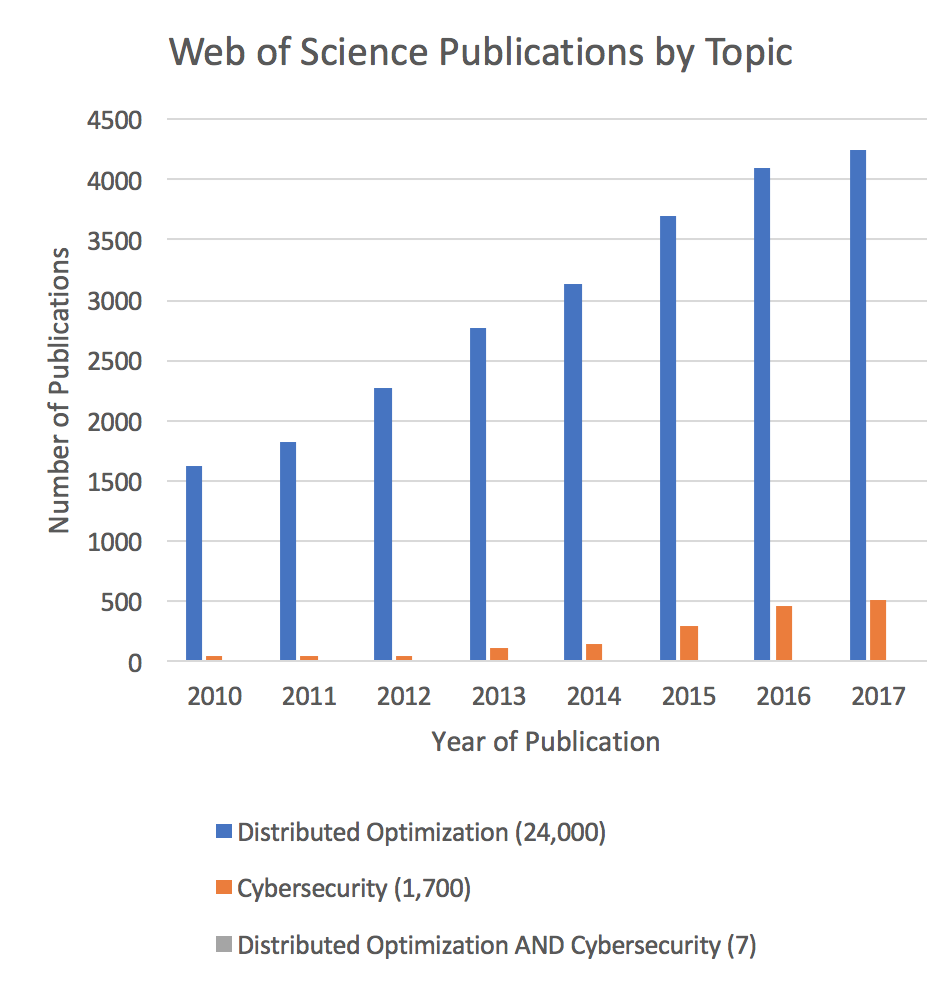}
	\caption{Publication frequency for the topics of `Distributed Optimization', `Cybersecurity', and the intersection of the two fields, 2010-2017. Legend captions include the total number of publications over this period; the number of publications at the intersection is three orders of magnitude less than each field.}
	\label{fig:citations}
\end{figure}

Prompted by these concerns, cybersecurity research has increased in recent years, as shown in Figure \ref{fig:citations}. However, this research is far outpaced by the development of new distributed optimization algorithms, and research on the intersection of the two fields - the security of distributed optimization algorithms- is still nascent.  Instead, the most cutting-edge research on distributed optimization techniques under adversarial attacks is found not in power systems or optimal control, but in the field of artificial intelligence.  

The deployment of distributed machine learning algorithms to very large computation networks or \textit{federated machine learning} platforms (in which consumer feedback e.g. in smartphone apps directly trains an estimator) has led researchers to work on algorithms which can be resilient to software crashes, network dropouts, and adversarial attacks on local nodes. 

Most of these machine learning problems are solved with consensus algorithms, in which local nodes solve private problems and reach consensus on a global estimator through iteration with a fusion node (aggregator-based model) or with their neighbors (fully-decentralized model). Although widely used, it has been shown that arbitrary attacks by a single corrupted node can lead consensus algorithms to fail \cite{Sundaram2016, BlanchardEPFL2017}- making it critical to understand potential attack vectors, detection algorithms, and mitigation strategies. 

As a result, research to combat these approaches has considered a number of methods for developing techniques to allow the consensus algorithm to converge to a best estimate under attack \cite{Sundaram2016}:
 \begin{itemize}
 \item \textit{Filtering techniques} remove outliers from the set of proposed updates \cite{Zhang2014, Sundaram2016, Liao2018}, in extreme cases using just the median estimator \cite{Xie2018, Yin2018} in order to provide tolerance to an attack of half of the nodes. These approaches can be computationally intensive, and may require that the local update functions are all drawn from the same distribution- feasible in statistical estimation, but unlikely in device scheduling.
 \item \textit{Nonlinear weighting schemes} take advantage of all estimators, but dynamically scale down the impact of suspicious updates. While these are guaranteed to move towards optimum \cite{BlanchardEPFL2017, Su2016a, Chen2018}, they have a larger optimality gap than other techniques \cite{Alistarh2018}.
 \item  \textit{Round-robin techniques} seek to detect and remove attacked updates by iteratively computing the consensus estimate with different combinations of dropped-out nodes, and using the most stable estimate \cite{Nabavi2015, Liao2016, Liao2018}. This is algorithmically simple but computationally intensive, and is most feasible when the number of attackers is well-bounded.
 \end{itemize}
 
 Additionally, \cite{Alistarh2018} demonstrates the improved convergence of these algorithms when compromised nodes can be switched off. 
 
 While many of these algorithms were developed for aggregator-coordinated machine learning problems, many of them can be adapted to fully-decentralized structures, as in the estimation problems studied in \cite{Chen2017b, Chen2018} and the optimization problem studied in  \cite{Nabavi2015}. 

Some of these techniques are beginning to also be deployed in smartgrid optimization problems: \cite{Weldehawaryat2017} demonstrates the hurdle which compromised nodes create for scheduling, much like \cite{Sundaram2016a} does for generalized consensus problems.

However, most research on cybersecurity in power system applications has not considered scheduling and optimization problems, but rather state estimation in transmission networks. Liu et al demonstrate in \cite{Liu2011}  that a single compromised node could introduce an arbitrarily large estimation error, similar to the generalized results in \cite{Sundaram2016a}. Subsequent research demonstrated the limitations of attack detection and mitigation strategies under a number of centralized state estimation models \cite{Li2012, Ozay2013, Liang2017}. 

 These approaches have also been extended to distributed optimization techniques in \cite{Liao2016}, where a round-robin ADMM method is used to identify compromised nodes; this approach is demonstrated experimentally in \cite{JanLiu2018}.   Similar research has also developed fully-decentralized algorithms for state estimation, with \cite{Kekatos2013} utilizing compressed sensing techniques to recover an estimate of system state when noise injection is bounded, though \cite{Vukovic2014} demonstrates that these approaches still fail when a node is able to inject arbitrary noise.

\subsection*{Novel Contributions}

This work advances prior literature by providing an attack detection algorithm for the \textit{alternating direction method of multipliers}, a popular optimization algorithm used for convex optimization and consensus problems. Only requiring that the private objective functions be convex, the attack detection algorithm works for both constrained and unconstrained problems, and both aggregator-coordinated and fully-decentralized architectures. By directly identifying compromised nodes, we bypass many of the objections to the filtering, nonlinear averaging, and round-robin techniques described above. Additionally, the proposed detection algorithm can readily be integrated into computational techniques, such as those described in \cite{Alistarh2018}.

In developing this, this paper provides the following novel contributions to existing literature:
\begin{itemize}
	\item Outline a taxonomy of attack vectors for decentralized optimization algorithms
	\item Detail attack detection, localization, and mitigation techinques for the alternating direction method of multipliers
	\item Demonstrate and verify an algorithm for detecting noise-injection attacks in the alternating direction method of multipliers
	\item Outline the unique security challenges of fully-decentralized optimization
	\item Describe potential architectures for security in fully-decentralized optimization
\end{itemize}

\section{Outline}

The remainder of this article is organized as follows. In Section IV we establish some mathematical background and notation. Section V outlines a taxonomy of methods by which an attacker may seek to compromise a decentralized optimization algorithm, and briefly discuss challenges to detection, localization, and mitigation.

Section VI derives an algorithm for detecting noise-injection attacks in convex optimization problems solved with ADMM, and Section VII provides results for a set of simulations with randomly-generated quadratic programs (QPs). Sections VIII and IX describe limitations and extensions, respectively.

Sections X and XI extend the cyberattack detection algorithm by considering the unique security challenges presented by fully-decentralized optimization algorithms. Section XII outlines architectures which can be used to overcome the weaknesses inherent in fully-decentralized architectures. Finally, Section XIII summarizes the article's main results.

\section{Mathematical Background and Notation}

\subsection{Notation}
Without loss of generality, we focus on an aggregator-coordinated system with two nodes which compute the $x$-update and $z$-update steps in a decentralized optimization problem. Also without loss of generality, we will consider an attack in which the $x$-update node has been compromised, and use the following notation:

\begin{itemize}
	\item $A,B$ Constraint matrices in ADMM binding constraint
	\item $c,d$ Linear cost vectors
	\item $f(x),g(z)$ generalized objective functions
	\item $H$ Hessian
	\item $i,j$ iterates
	\item $k$ iterate limit
	\item $m,n$ Number of dimensions of $z$ and $x$ variables respectively
	\item $p$ Number of binding constraints
	\item $P, Q$ quadratic cost functions in sample problems
	\item $r, s$ dimensions in Hessian
	\item $u$ Scaled dual variable
	\item $w$ combined variable for centralized solution
	\item $x,z$ Optimization variables
	\item $y$ Unscaled dual variable
	\item $\mathcal{X,Z}$ The private constraints, knowable only to the compute nodes and not publicly shared 
	\item $x^{\star k}$ The unattacked update, solving $x^{k} := \; \text{argmin}_{x \in \mathcal{X}} \; f(x) + \frac{\rho}{2}\lVert Ax + Bz^{k-1} - c + u_z^{k-1}\rVert^2_2$
	\item $\tilde{x}^k$ The attacked signal provided by the $z$ node
	\item $x^{? k}$ The variable update received by the $z$-update node, of unknown validity
	\item $\hat{x}^k$ A best response created by the $z$-update node which has received a signal perceived as being an attack
\end{itemize}

\subsection{Decentralized Optimization}

In general, decentralized optimization problems can be cast as:
\begin{align*}
x^*, y^* = \text{argmin}_{x,z} & \quad W\left(  x, z \right)\\
\text{s. to:}&\quad x, z \in \mathcal{W}\\ 
\end{align*}
where $x$ and $z$ further satisfy local problems:
\begin{align*}
x^* = \text{argmin}_{x} &\quad f(x, z^*) \\
&\quad  x \in \mathcal{X}\\
z^* = \text{argmin}_{x} &\quad g(z, x^*) \\
&\quad  z \in \mathcal{Z}\\
\end{align*}

In this discussion we focus on iterative methods, where updates $x^k = \text{argmin}_{x\in\mathcal{X}} \, f(x, z^{k-1})$, $z^k =\text{argmin}_{z\in\mathcal{Z}} \, g(z, x^{k-1})$ solve local updates in an algorithm which converges to a global solution, i.e. $x^k \rightarrow x^*, z^k \rightarrow z^*$ as $k\rightarrow \infty$.

In addition to computational benefits, this architecture is also advantageous when the local objective functions $f(x),\,g(z)$ or local constraint sets $x \in \mathcal{X}, \, z \in \mathcal{Z}$ contain \textit{private} information which can not be directly shared with other participants or the aggregator. In this scenario, the updates $x^k, z^k$ do not directly reveal private information, yet still allow the system to converge to the global optimum.  This is common in markets and scheduling problems, where participants have private constraints or utility functions which they do not wish to share for economic or security reasons.

Additionally, we assume that some information about the private constraint set is publicly knowable, creating a superset which contains the private constraint set: $\mathcal{X} \subset \mathcal{X}_{\text{pub}},\, \mathcal{Z} \subset \mathcal{Z}_{\text{pub}}$. Note that as all variable updates are in the private constraint set, they must also be in the public set: $x^k \in \mathcal{X} \in \mathcal{X}_{\text{pub}}$.

\subsection{The Alternating Direction Method of Multipliers}

Although any distributed optimization algorithm may be subject to cyberattack, we specifically consider the \textit{Alternating Direction Method of Multipliers} (or ADMM) algorithm reviewed in \cite{Boyd2010}, which has gained popularity due to its simple formulation and guarantees of convergence for convex problems. The ADMM algorithm is used to decompose a problem with separable objective and constraints:

\begin{align}
\min_{x,z} &\quad f(x) + g(z)\\
\text{s.t.}&\quad Ax+Bz = c \label{eqn:bindingconstraint}\\ 
&\quad  x \in \mathcal{X}\\
&\quad  z \in \mathcal{Z}
\end{align}
Note that only equation \ref{eqn:bindingconstraint} links $x$ and $z$. 

Adding a penalty term for violations of the linking constraint, we can create an augmented Lagrangian defined in the domain $\{ x\in\mathcal{X}, z\in\mathcal{Z}\}$:
\begin{align}
&L_{\rho}(x,z,y) =  \\
&f(x)+g(z)+y^T(Ax+Bz-c) + (\rho/2)||Ax+Bz-c||_2^2 \nonumber
\end{align}
The standard ADMM algorithm can then be expressed with respect to this augmented Lagrangean $L_{\rho}(\cdot)$ as:
\begin{align*}
x^{k+1} &:= \quad \text{argmin}_{x \in \mathcal{X}} \; L_{\rho}(x      ,z^k,y^k) \\
z^{k+1} &:= \quad \text{argmin}_{z \in \mathcal{Z}} \; L_{\rho}(x^{k+1},z  ,y^k) \\
y^{k+1} &:= \quad y^k + \rho(Ax^{k+1} + Bz^{k+1} - c)
\end{align*}
or in scaled form:
\begin{align}
x^{k+1} &:= \quad \text{argmin}_{x \in \mathcal{X}} \; f(x) + \frac{\rho}{2}\lVert Ax + Bz^k - c + u^k \rVert^2_2 \\
z^{k+1} &:= \quad \text{argmin}_{z \in \mathcal{Z}} \; g(z) + \frac{\rho}{2}\lVert Ax^{k+1} + Bz - c + u^k \rVert^2_2 \\
u^{k+1} &:= \quad u^k + Ax^{k+1} + Bz^{k+1} - c
\end{align}

The iterations stop when:
\begin{itemize}
	\item The primal residual $r^k = Ax^k + Bz^k - c$ has a magnitude below a threshold $\epsilon_{\text{pri}}$, i.e. $||r^k||_2 \leq \epsilon_{\text{pri}}$
	\item The dual residual $s^k = \rho A^T B (z^k - z^{k-1})$ has a magnitude below a threshold $\epsilon_{\text{dual}}$, i.e. $||s^k||_2 \leq \epsilon_{\text{dual}}$
\end{itemize}

As described in \cite{Boyd2010} the ADMM iterations will converge to the optimal value of the objective function, and the primal and dual residuals will converge to zero.

In this ADMM formulation, we refer to the ``aggregator'' (or ``fusion node'') as the agent responsible for updating $u$. The aggregator:
\begin{enumerate}
	\item receives updates from the $x$-update and $z$-update steps, 
	\item computes the $u$-update
	\item broadcasts the updated value of $u$ to the nodes responsible for computing the $x$- and $z$-updates.
\end{enumerate}

\subsubsection{Consensus problems}
ADMM has become a popular tool for solving \textit{consensus problems} in both machine learning and power systems engineering, as both fields deal with problems where computation may be spread across thousands (or millions) of nodes.  In a consensus problem, local variables $x_i$ and objective functions $f_i(x_i), i\in \{1,\ldots,N\}$ are united by the constraint that at optimality, the local variables mirror the global variable $z$:
\begin{align*}
\min_{x,z} &\quad \sum_{i=1}^N f_i(x_i) \\
\text{s.t.}&\quad x_i = z, \quad i = 1, \ldots,N 
\end{align*}
Collected, the ADMM form of this is:
\begin{align*}
x^{k+1}_i &:= \quad \text{argmin}_{x_i \in \mathcal{X_i}} \; f_i(x_i)+y_i^{kT}(x_i-z^k) + \frac{\rho}{2} ||x_i - z^k||^2_2 \\
z^{k+1} &:= \quad \frac{1}{N} \sum_{i=1}^N x_i^{k+1} + (1/\rho) y_i^k\\
y^{k+1}_i &:= \quad y_i^k + \rho(x_i^{k+1} - z^{k+1} )
\end{align*}

Note that this structure generalizes minibatch stochastic gradient descent, in which the local objectives are the gradients of the local estimators.

As this consensus algorithm is a specific example of ADMM, results for the general ADMM problem will also apply to the ADMM consensus algorithm. For clarity, in this paper we will continue to refer to $x$- and $z$-update nodes, even though a consensus problem may have thousands or millions of $x$-nodes and a single $z$-update (aggregator) node. Despite this difference in scale, the results we derive for the general 2-node system will still be applicable to the consensus problem.

We will also note that consensus problems are of particular interest in many power system and machine learning problems, as well as in fully-decentralized optimization problems, which include consensus networks. 

\section{Attack Vectors in Decentralized Optimization}\label{sec:attackvector_overview}
When local problems contain private information, the central coordinator can only check $x,z \in \mathcal{W}$ and cannot directly verify that the updates $x^{?k}, z^{?k}$ in fact solve the private optimization problems. In this scenario, a malicious node may submit a distorted update $\tilde{x}^k$ in order to mislead the central coordinator with the goal of creating a sub-optimal solution, an infeasible solution, or prevent convergence of the iterative algorithm.

In this section, we briefly outline potential attack vectors and methods for addressing them. Of these attacks, the most generalizable but difficult to identify is zero-mean noise injection, which we consider in detail below. These attack vectors can be combined, but detection and mitigation efforts will usually address these separately.

\subsubsection*{Sub-optimal solution} 
The compromised node solves a modified objective function $\tilde{f}(x)$ resulting in $f(x^k) < \tilde{f}(x^k)$ and consequently $W(\tilde{x}^k,z^k) > W(x^k,z^k)$ as $k \rightarrow \infty$. As this does not change the problem structure (e.g. convexity, private constraints), it is difficult to discern from an equivalent problem with slightly modified characteristics, e.g. different consumer preferences. This makes attack prevention a matter of appropriately structuring game-theoretic incentives to avoid malicious distortion. As such, we do not consider this in greater detail here. 
	
\subsubsection*{Infeasible Private Constraints} 
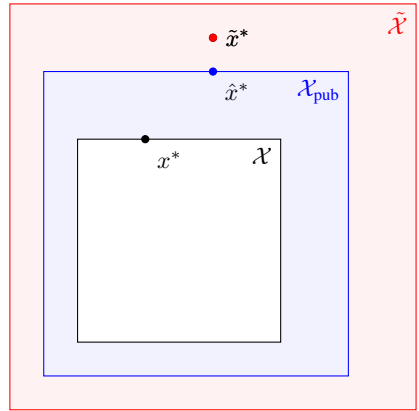
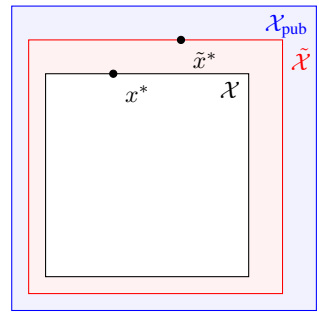
\begin{figure}
	\subfloat[Attacked system in which detection and (limited) mitigation is possible, as $\mathcal{X} \subset \mathcal{X}_{\text{pub}} \subset \tilde{\mathcal{X}}$ .]{
	\begin{tikzpicture}[
	scale=0.9, 
	every node/.style={transform shape},
	bluenode/.style={circle, draw=blue!60, fill=blue!5, very thick, minimum size=7mm},
	rednode/.style={circle, draw=red!60, fill=red!5, very thick, minimum size=5mm},
	]
	\draw[red,fill=red!5] (0,0) rectangle (6,6) node[anchor=north east] {$\tilde{\mathcal{X}}$};
	\draw[blue,fill=blue!5] (0.5,0.5) rectangle (5,5) node[anchor=north east] {$\mathcal{X}_{\text{pub}}$};
	\draw[black,fill=white] (1,1) rectangle (4,4) node[anchor=north east] {$\mathcal{X}$};
	
	\filldraw (2,4) circle[radius=1.5pt];
	\node[below right=1pt of {(2,4)}, outer sep=1pt] {$x^*$};
	\filldraw (3.,5.5) circle[radius=1.5pt];
	\node[right=1pt of {(3.,5.5)}, outer sep=1pt] {$\tilde{x}^*$};
	\filldraw[red] (3.,5.5) circle[radius=1.5pt];
	\node[right=1pt of {(3.,5.5)}, outer sep=1pt] {$\tilde{x}^*$};
	\filldraw[blue] (3.0,5.0) circle[radius=1.5pt];
	\node[below right=1pt of {(3.,5)}, outer sep=1pt] {$\hat{x}^*$};
	\end{tikzpicture}
	\label{fig:attackedDetectFeasible}}
	\qquad
	\subfloat[Attacked system in which detection and mitigation is not possible, as the distorted constraint set is a subset of the publicly known bounds on the constraints,  $\mathcal{X} \subset \tilde{\mathcal{X}} \subset  \mathcal{X}_{\text{pub}} $.]{
	\begin{tikzpicture}[
	scale=0.9, 
	every node/.style={transform shape},
	bluenode/.style={circle, draw=blue!60, fill=blue!5, very thick, minimum size=7mm},
	rednode/.style={circle, draw=red!60, fill=red!5, very thick, minimum size=5mm},
	]
	\draw[white,fill=white] (0,0) rectangle (6,6) node[anchor=north east] {$\mathcal{X}_{\text{pub}}$};
	\draw[blue,fill=blue!5] (0.5,0.5) rectangle (5,5) node[anchor=north east] {$\mathcal{X}_{\text{pub}}$};
	\draw[red,fill=red!5] (0.75,0.75) rectangle (4.5,4.5) node[anchor=north west] {$\tilde{\mathcal{X}}$};
	\draw[black,fill=white] (1,1) rectangle (4,4) node[anchor=north east] {$\mathcal{X}$};
	
	\filldraw (2,4) circle[radius=1.5pt];
	\node[below right=1pt of {(2,4)}, outer sep=1pt] {$x^*$};
	\filldraw (3.,4.5) circle[radius=1.5pt];
	\node[below right=1pt of {(3.,4.5)}, outer sep=1pt] {$\tilde{x}^*$};
	
	\end{tikzpicture}
	\label{fig:attackedDetectInfeasible}}
	\caption{Graphical example of how an attacker may distort the constraint set from $\mathcal{X}$ to $\tilde{\mathcal{X}}$ to create an optimum outside of the truly feasible set, and the limited ability to mitigate these impacts by projecting onto a publicly knowable constraint set $\mathcal{X}_{\text{pub}}$.}
	\label{fig:infeasibleconstraints}
\end{figure}

An attacked node may replace the constraints $x\in \mathcal{X}$ with a modified constraint set $x\in \tilde{\mathcal{X}}$ in order to result in an update which is not feasible $\tilde{x}^k \not\in \mathcal{X}$, as shown in Figure \ref{fig:infeasibleconstraints}. This leads the system to converge to a point which does not reflect actual conditions, e.g. a schedule which is not operationally feasible or a machine learning estimator distorted by false data.
Because these constraints may arise from stochastic processes (user preferences, input data) the aggregator cannot directly discern an attack from an unattacked update with different underlying data. As a result, defenses from this attack are limited to cases where the aggregator can bound the support of the problem, i.e. a publicly knowable constraint set $\mathcal{X}_\text{pub}$. In this scenario, attacks can be detected when updates lie out of the support region, and the effect of the attack can be mitigated by projecting the  update back onto the support region.

\subsubsection*{Infeasible Linking Constraint}
Knowing that a solution to the master problem must satisfy $Ax + Bz = c$, an attacker may present an update which it knows to be unreachable with these constraints, i.e. $ \nexists z \in \{\mathcal{Z} | A\tilde{x}^k + Bz = c\}$. In this case, convergence stops as the primal residual  $r^k = ||Ax^k + Bz^k - c||^2_2$ remains nonzero. 
This attack relies on th attacker knowing some bounds on the set reachable by its neighbor $\mathcal{Z} \subset \mathcal{Z}_{\text{pub}} $ in order to shape this attack. Similar to above, detection relies on the same knowledge of the public constraint set, and mitigation relies on projecting the attacked signal back onto the feasible space. 
If the aggregator has knowledge of the public constraints $\mathcal{Z}_{\text{pub}}$, then detection and mitigation can easily be implemented at the aggregator level for each signal.

\subsubsection*{Non-Convergence} 
\begin{figure}
	\centering
	\includegraphics[scale = 0.5]{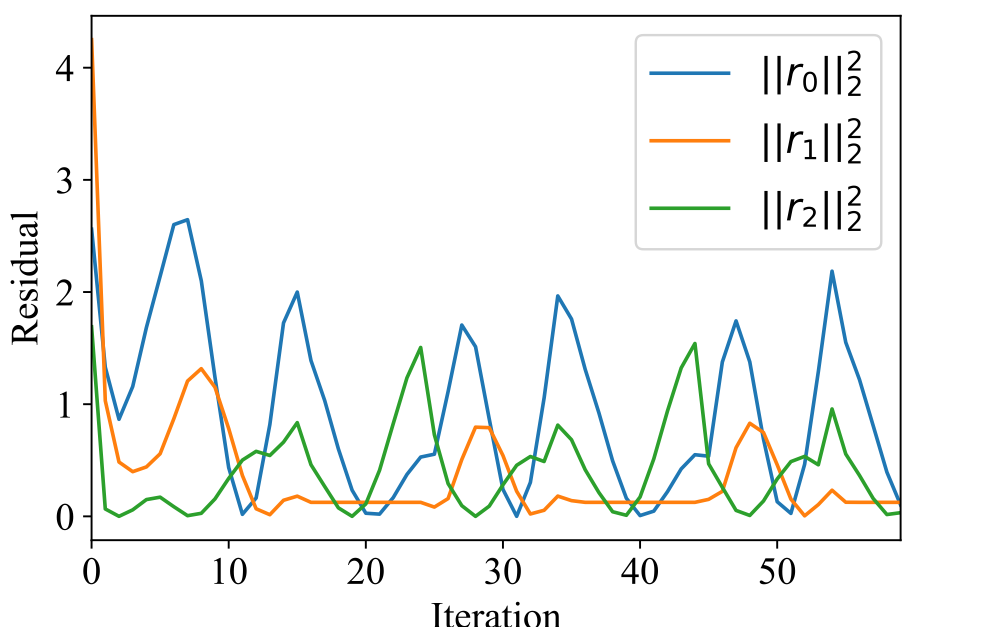}
	\caption{Plot of local residuals in a 3-node system under noise-injection attack. Injection of noise prevents problem convergence at nodes across the network.}
	\label{fig:noiseinjectionresiduals}
\end{figure}

An attacker may add a noise term which varies with each iteration, i.e. $\tilde{x}^i = x^{\star i} + \delta(i)$. When the noise term is larger than the rate of convergence of the problem, this results in a persistent residual which remains above the convergence tolerance, as shown in Figure \ref{fig:noiseinjectionresiduals}.
We examine in detail the case where the noise is unbiased (zero-mean). When the noise is biased, the problem can be considered a combination of a zero-mean noise injection attack and one of the attacks described above, and addressed accordingly.

\section{Developing a Detection Algorithm for Noise-Injection Attacks}\label{sec:algorithm_noisedetection}

Because the ADMM algorithm relies on private actors computing $x-$ and $z-$updates, these actors can distort the prolem by providing inaccurate updates, with the goal of creating a suboptimal or infeasible solution or preventing convergence. We study the case of distorting the $x$-update through injection of zero-mean noise, which is intended to prevent convergence of the algorithm.  The new update becomes:

$$
x^{i+1} := \quad \text{argmin}_{x \in \mathcal{X}} \quad f(x) + \frac{\rho}{2}\lVert Ax + Bz^i - c + u^i \rVert^2_2  + \delta(i)
$$
Where $\delta(i)$ is a noise term which changes on each iteration. The $z$-actor or the system aggregator is challenged to identify the attack and take preventative actions before convergence is prevented.

\subsection{Attack Detection Overview}

We play the part of the aggregator or z-actor, and wish to detect an attack of $x$, using the values of the $x^i, z^i, u^i,\; i=1\ldots k$ iterates. We rely on \textit{a priori} knowledge that $f(x)$ must be convex, we trust the computations of $z^i$, and we can verify the $u^i$ values directly by using the other iterates.

We will detect attacks by assessing the convexity of $f(x)$ implied by the $x^i$ iterates, \textit{without being able to directly assess $f(x)$}. 

\begin{figure}
	\centering
	\subfloat[]{
		\includegraphics[scale = 0.32]{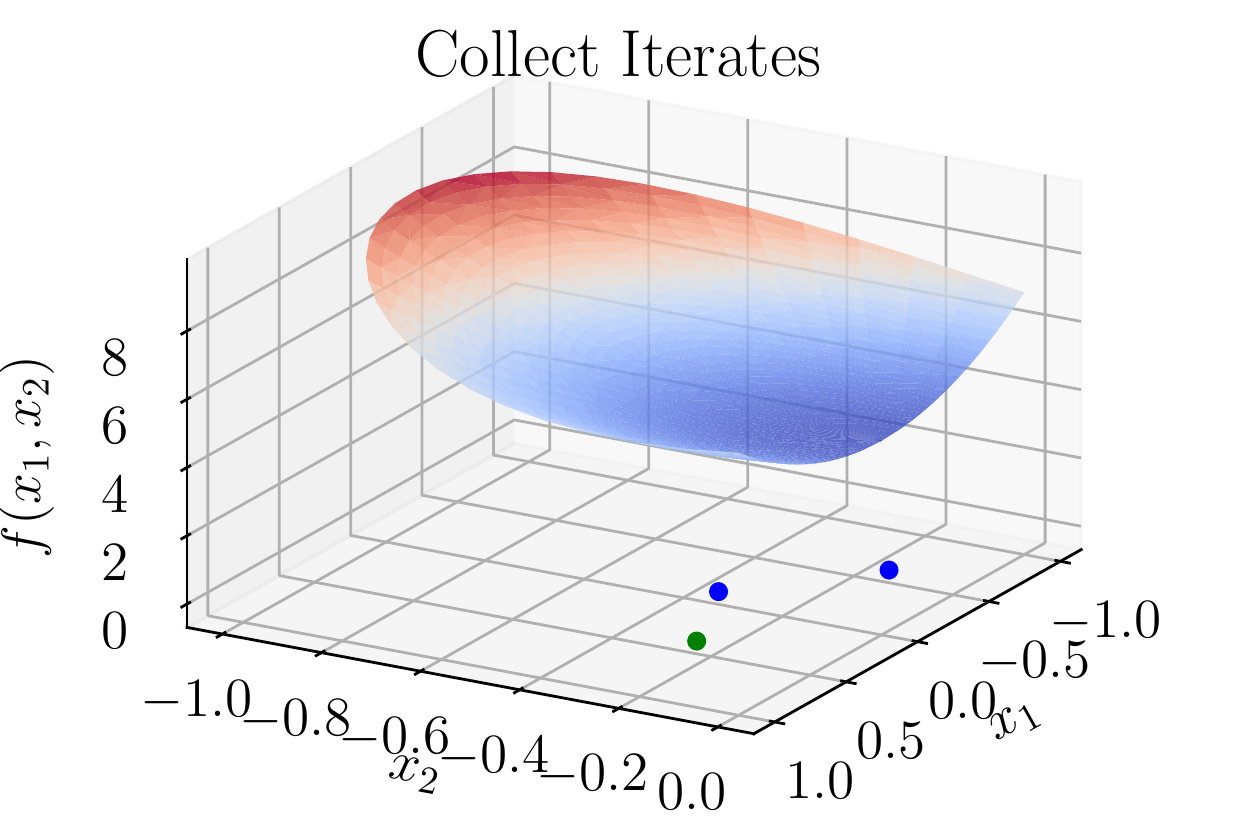}
		\label{fig:conceptModelIterates}
	}
	\subfloat[]{
		\includegraphics[scale = 0.32]{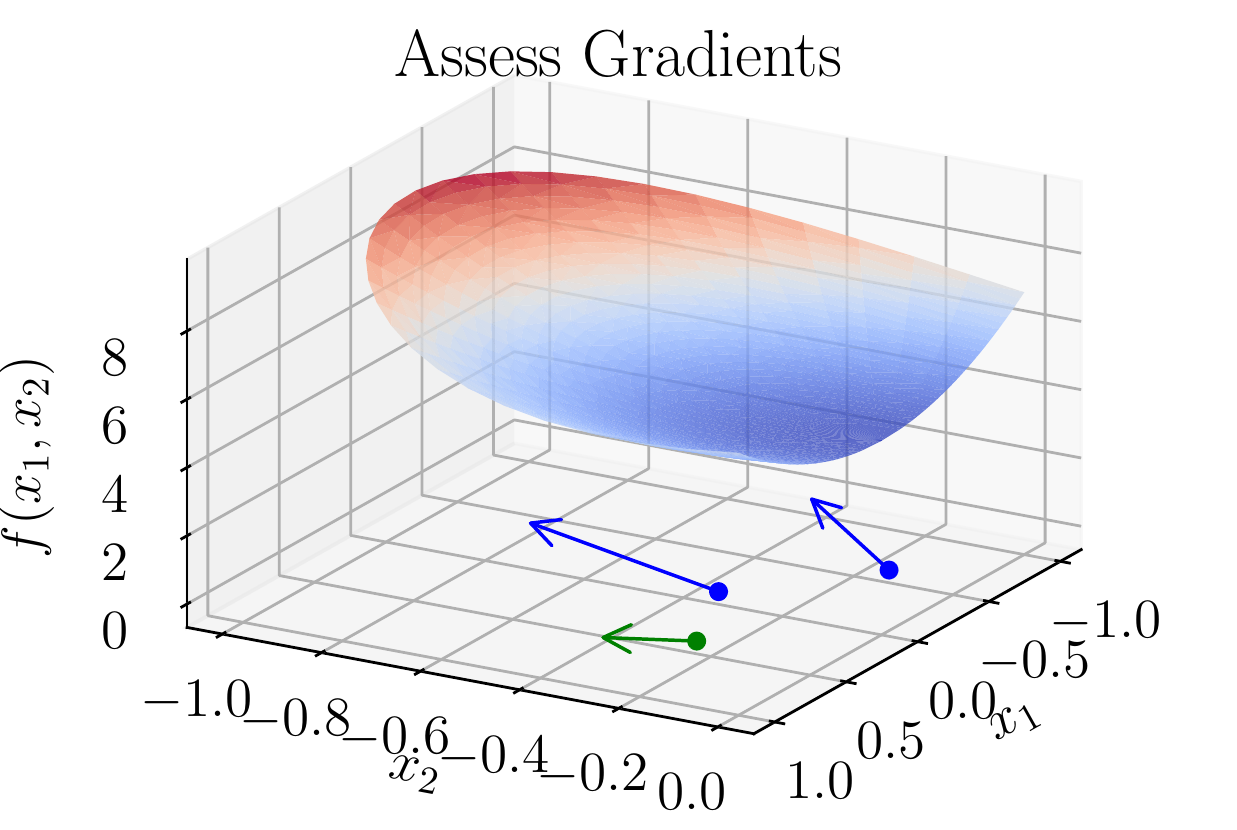}
		\label{fig:conceptModelGradients}
	}
	\ifarticle
	
	\fi
	\subfloat[]{
		\includegraphics[scale = 0.32]{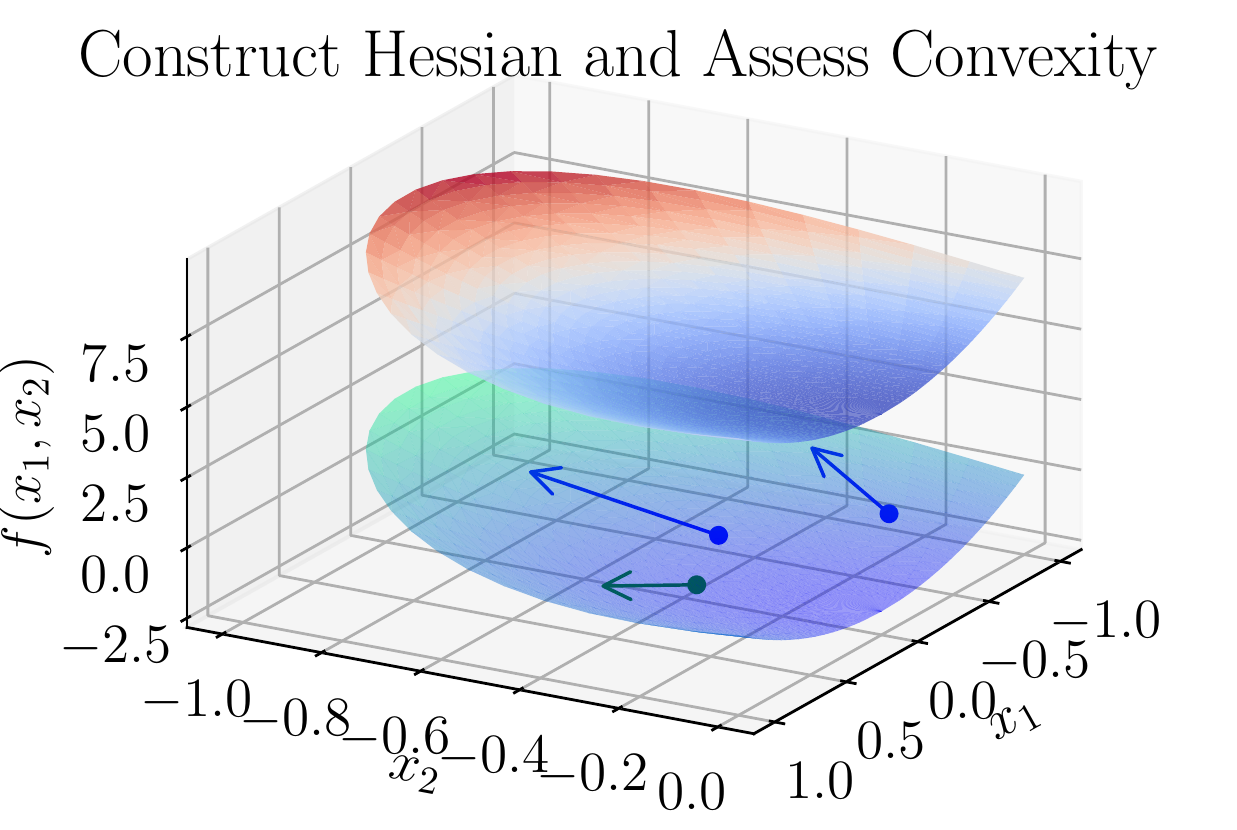}
		\label{fig:conceptModelConvexity}
	}
	\caption{Conceptualization of the attack detection process for a $n=2$ toy problem. Although the validator cannot directly assess the private objective function $f(x_1,x_2)$ (shown in red-blue gradient), they are provided with a sequence of iterates $x^k$, and can  to use a subset of these for the detection algorithm, shown here as $(x_1,x_2)$ points with $f(x)=0$ due to the unknown objective value.  Using information from the sequence of iterates, the validator can then assess the implied gradient of $f(x)$ at each of these points. Finally, from these gradients, the validator can estimate the Hessian, and use this to evaluate the convexity of the (unknown) objective function. This can be conceptualized as a local quadratic approximation of the objective function at the reference point, shown in blue-green.}
	\label{fig:unattackedVectors}
\end{figure}

To do this, we will use the $z-$ and $u$-update values to evaluate the local gradient of $f(x)$, then use these local gradients to construct a finite-differences approximation of the Hessian of $f(x)$. Testing whether this Hessian is symmetric positive semi-definite will then allow us to test the convexity of the $x$-updates; this is visualized in Figure \ref{fig:unattackedVectors}. Given our \textit{a priori} knowledge that $f(x)$ must be convex, any updates which result in a nonconvex Hessian must represent an attack.

The algorithm progresses in three steps:
\begin{itemize}
	\item Assess gradient
	\item Construct Hessian
	\item  Evaluate eigenvalues of Hessian
\end{itemize}

\subsubsection{Assess Gradient}

We rely on the proof of ADMM convergence\footnote{specifically building on the \textit{Proof of Inequality A2} on pg 108 of \cite{Boyd2010}} provided by \cite{Boyd2010}, and while we focus on the x-update the same approach can be used for security checks conducted by the $x$-update node.  The gradient of $f(x)$ evaluated at $x^i$ is found as:
\begin{align*}
0 &\in \partial L_{\rho}(x^i, z^{i-1}, u^{i-1}) \\
0 &\in \partial f(x^i) + A^T + \rho A^T(A x^i + Bz^{i-1} - c)\\
0 &\in \partial f(x^i) + A^T (y^i - \rho B(z^i-z^{i-1}))\\
\partial f(x^i) &= -A^T (y^i - \rho B(z^i-z^{i-1}))\\
\partial f(x^i) &= -A^T (\rho u^i - \rho B(z^i-z^{i-1}))\\
\end{align*}
where the last equality simply uses the scaled dual variable. 

We will interchangeably use the notation $f(x^i)$ and $f(x)|^{x^i}$ to both indicate the definite evaluation of $f(x)$ at the point $x^i$. The iterates $i=1,\ldots,k$ thus give us a sequence of gradient evaluations at points $x^i, i=1,\ldots,k$. From this sequence of gradients, we wish to construct the Hessian.

\subsubsection{Construct Hessian}

We wish to construct the Hessian matrix
$$
H = 
\begin{bmatrix}
\dfrac{\partial}{\partial^2 x_1}f(x) & \dfrac{\partial}{\partial x_1 \partial x_2}f(x) & \ldots & \dfrac{\partial }{\partial x_1 \partial x_n}f(x) \\
\dfrac{\partial}{\partial x_2 \partial x_1}f(x) & \dfrac{\partial}{\partial^2 x_2}f(x) & \ldots & \dfrac{\partial}{\partial x_2 \partial x_n}f(x) \\
\vdots  & \vdots & \ddots & \ldots \\
\dfrac{\partial}{\partial x_n \partial x_1}f(x) & \dfrac{\partial}{\partial x_n \partial x_2}f(x) & \ldots & \dfrac{\partial}{\partial^2 x_n}f(x) \\
\end{bmatrix} 
$$

To construct the Hessian, we utilize information gained from the gradient evaluations at each of the iterates to compute a numeric approximation of the local Hessian, using a Taylor series expansion of the Hessian around that point. 

\subsubsection*{Example: 2-Dimension Case}
\begin{figure}
	\centering
	\includegraphics[scale = 0.7]{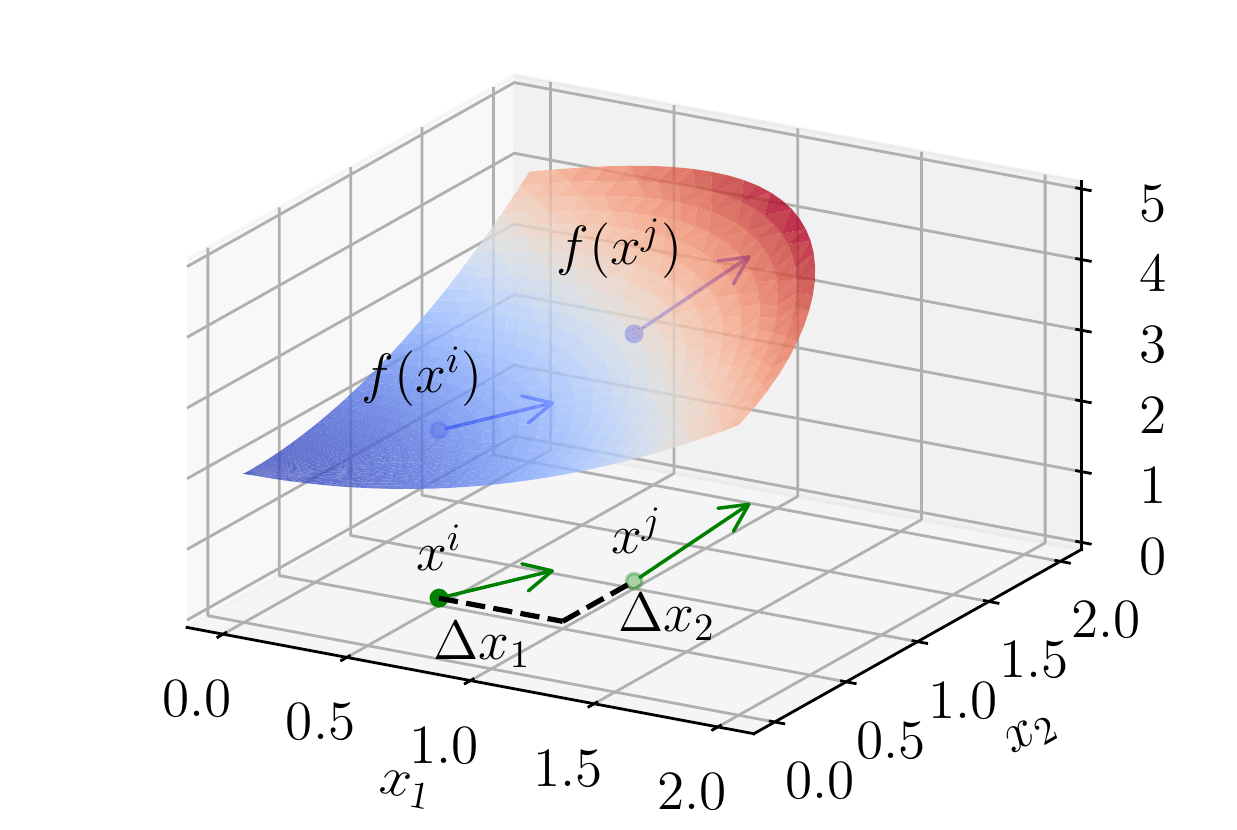}
	\caption{Visualization of the relationship between the iterates $x^i, x^j$, the function surface $f(x)$, and the gradient evaluations $\partial f(x)|^i, \partial f(x)|^j$. While the function surface $f(x)$ and its Hessian cannot be directly evaluated, the change in gradient evaluations between $x^i$ and $x^j$ can be used to derive an approximation of the Hessian.}
	\label{fig:gradientChanges}
\end{figure}

For clarity, we begin with a two-dimensional case, i.e. $x^i = [x_1^i, x_2^i]^T$ and visualized in Figure \ref{fig:gradientChanges}. We will consider two points $x^i$ and $x^j$, and will note the dimensions of these points using the subscripts $1$ and $2$.

In this case, the change in gradient from point $x^i$ to $x^j$ can be approximated using the following difference equations: 

\ifarticle 
\noindent\resizebox{\columnwidth}{!}{\parbox{\columnwidth}{
\fi
\begin{align*}
\frac{\partial}{\partial x_1}f(x)|^{x^j} &= \frac{\partial}{\partial x_1}f(x)|^{x^i} +  \frac{\partial}{\partial^2 x_1}f(x)|^{x^i}(x^j_1 - x^i_1) + \frac{\partial}{\partial x_1 \partial x_2}f(x)|^{x^i}(x^j_2 - x^i_2)\\
\frac{\partial}{\partial x_2}f(x)|^{x^j} &= \frac{\partial}{\partial x_2}f(x)|^{x^i} +  \frac{\partial}{\partial^2 x_2}f(x)|^{x^i}(x^j_2 - x^i_2) + \frac{\partial}{\partial x_2 \partial x_1}f(x)|^{x^i}(x^j_1 - x^i_1)\\
\end{align*}
\ifarticle
}}
\fi

This is visualized in Figure \ref{fig:gradientChanges}, where $\Delta x_1 = (x^j_1 - x^i_1)$ and $\Delta x_2 = (x^j_2 - x^i_2)$. From the previous step, we can directly evaluate the gradient terms 
\begin{align}
\frac{\partial}{\partial x_r}f(x)|^{x^k}, \quad & r\in{1,2}, k\in{i,j}
\end{align}
and wish to solve for the second-order terms, which will become the entries in the Hessian matrix. However, we have $n^2$ unknown second-order terms but only $n$ equations, making this system underdefined. In order to have a fully-defined system, we need another set of difference equations, found by considering the difference with respect to another point $k$. Using the notation $(x_1)|^k_j = x_1^k - x_1^j$ and $\frac{\partial}{\partial x_1}f(x)|^k_j = \frac{\partial}{\partial x_1}f(x)|^k - \frac{\partial}{\partial x_1}f(x)|^j$ we can now rearrange the problem into the matrix equation
\begin{align}
\overrightarrow{G} &= D \overrightarrow{H}
\end{align}
or for the 2-D case:
\begin{equation*}
\ifarticle 
\resizebox{\columnwidth}{!}{$
\fi
\begin{bmatrix}
\dfrac{\partial}{\partial x_1}f(x)|^k_j \\
\dfrac{\partial}{\partial x_2}f(x)|^k_j \\
\dfrac{\partial}{\partial x_1}f(x)|^k_i \\
\dfrac{\partial}{\partial x_2}f(x)|^k_i \end{bmatrix} = 
\begin{bmatrix}
(x_1)|^k_j & (x_2)|^k_j & 0 & 0\\
0 & 0 & (x_1)|^k_j & (x_2)|^k_j\\
(x_1)|^k_i & (x_2)|^k_i & 0 & 0\\
0 & 0 & (x_1)|^k_i & (x_2)|^k_i\\
\end{bmatrix}
\begin{bmatrix}
\dfrac{\partial}{\partial^2 x_1}f(x)\\
\dfrac{\partial}{\partial x_1 \partial x_2}f(x) \\
\dfrac{\partial}{\partial x_2 \partial x_1}f(x) \\
\dfrac{\partial}{\partial^2 x_2}f(x)\\
\end{bmatrix}%
\ifarticle
$}
\fi
\end{equation*}

In this form, $\overrightarrow{G}$ collects the changes in gradient evaluations induced by each difference equation, the matrix $D$ arranges the difference in coordinates at the evaluated iterates, and $\overrightarrow{H}$ unstacks the entries of the Hessian matrix $H$.  

The entries of the Hessian matrix can then be found as $ \overrightarrow{H} = D^{-1}\overrightarrow{G}$, or for the 2-D case:
\begin{equation*}
\ifarticle 
\resizebox{\columnwidth}{!}{$
\fi
\begin{bmatrix}
\dfrac{\partial}{\partial^2 x_1}f(x)\\
\dfrac{\partial}{\partial x_1 \partial x_2}f(x) \\
\dfrac{\partial}{\partial x_2 \partial x_1}f(x) \\
\dfrac{\partial}{\partial^2 x_2}f(x)\\
\end{bmatrix}
= 
\begin{bmatrix}
(x_1)|^k_j & (x_2)|^k_j & 0 & 0\\
0 & 0 & (x_1)|^k_j & (x_2)|^k_j\\
(x_1)|^k_i & (x_2)|^k_i & 0 & 0\\
0 & 0 & (x_1)|^k_i & (x_2)|^k_i\\
\end{bmatrix}^{-1}
\begin{bmatrix}
\dfrac{\partial}{\partial x_1}f(x)|^k_j \\
\dfrac{\partial}{\partial x_2}f(x)|^k_j \\
\dfrac{\partial}{\partial x_1}f(x)|^k_i \\
\dfrac{\partial}{\partial x_2}f(x)|^k_i \end{bmatrix}
\ifarticle
$}
\fi
\end{equation*}

\subsubsection*{Expanding to $n$ dimensions}
Generalized to $r,s \in l$ dimensions (where $l \leq n$), this becomes:
$$
\frac{\partial}{\partial x_r}f(x)|^{x^j} = \frac{\partial}{\partial x_r}f(x)|^{x^i} + \sum_{s=1}^l \frac{\partial}{\partial x_r \partial x_s}f(x)|^{x^i}(x^j_s - x^i_s)
$$

As there are $n$ dimensions in $x$, there are a total of $n^2$ unknown terms in the Hessian 
\footnote{or $n(n+1)/2$ if we utilize the fact that the Hessian is symmetric, i.e. that $\frac{\partial}{\partial x_1 \partial x_2} = \frac{\partial}{\partial x_2 \partial x_1}$. For simplicity, we will not take advantage of symmetry in this manuscript, though this symmetry can be utilized to accelerate the solution time of the algorithm.}. As each new point $x^j$ has $n$ dimensions, it adds a new set of $n$ difference equations. In order to fully define the set of equations for computing the Hessian, we need to collect a total of $n+1$ linearly independent points (reference point plus differences with $n$ points) to solve for the $n^2$ terms in $H$. Since computing the gradient utilizes values from two iterates, this means that for a problem with $n$ dimensions, the detection algorithm can be conducted on iterations $n+2$ and beyond.

In the following, we index dimensions by $1,2,\ldots,n$ and points by $a,b,\ldots,k$ where $k$ is used as the reference point. The vector $\overrightarrow{G}$ is thus composed by stacking the gradient evaluations:
$$
\renewcommand*{\arraystretch}{2.5}
\overrightarrow{G} = \begin{bmatrix}
\dfrac{\partial}{\partial x_1}f(x)|^k_a \\
\dfrac{\partial}{\partial x_2}f(x)|^k_a \\
\vdots \\ 
\dfrac{\partial}{\partial x_n}f(x)|^k_a \\
\vdots\\
\dfrac{\partial}{\partial x_1}f(x)|^k_{k-1} \\
\dfrac{\partial}{\partial x_2}f(x)|^k_{k-1}\\
\vdots \\ 
\dfrac{\partial}{\partial x_n}f(x)|^k_{k-1}\\
\end{bmatrix}
$$

The matrix $D$ is composed of a stack of $n$ block-diagonal matrices, each of which is the Kronecker product of the $n$-dimensional identity matrix and a set of coordinate differences from $x$-iterates. Each block of this is $n \times n^2$, producing  $D \in \mathbb{R}^{n^2 \times n^2}$:

$$
D = 
\begin{bmatrix}
I_n \otimes (x^k - x^a)^T \\
I_n \otimes (x^k - x^b)^T \\
\vdots \\
I_n \otimes (x^k - x^{k-1})^T \\
\end{bmatrix}
$$

We solve for the Hessian elements by computing $ \overrightarrow{H} = D^{-1}\overrightarrow{G}$ as before.

\begin{figure}
	\centering
	\subfloat[]{\includegraphics[scale = 0.5]{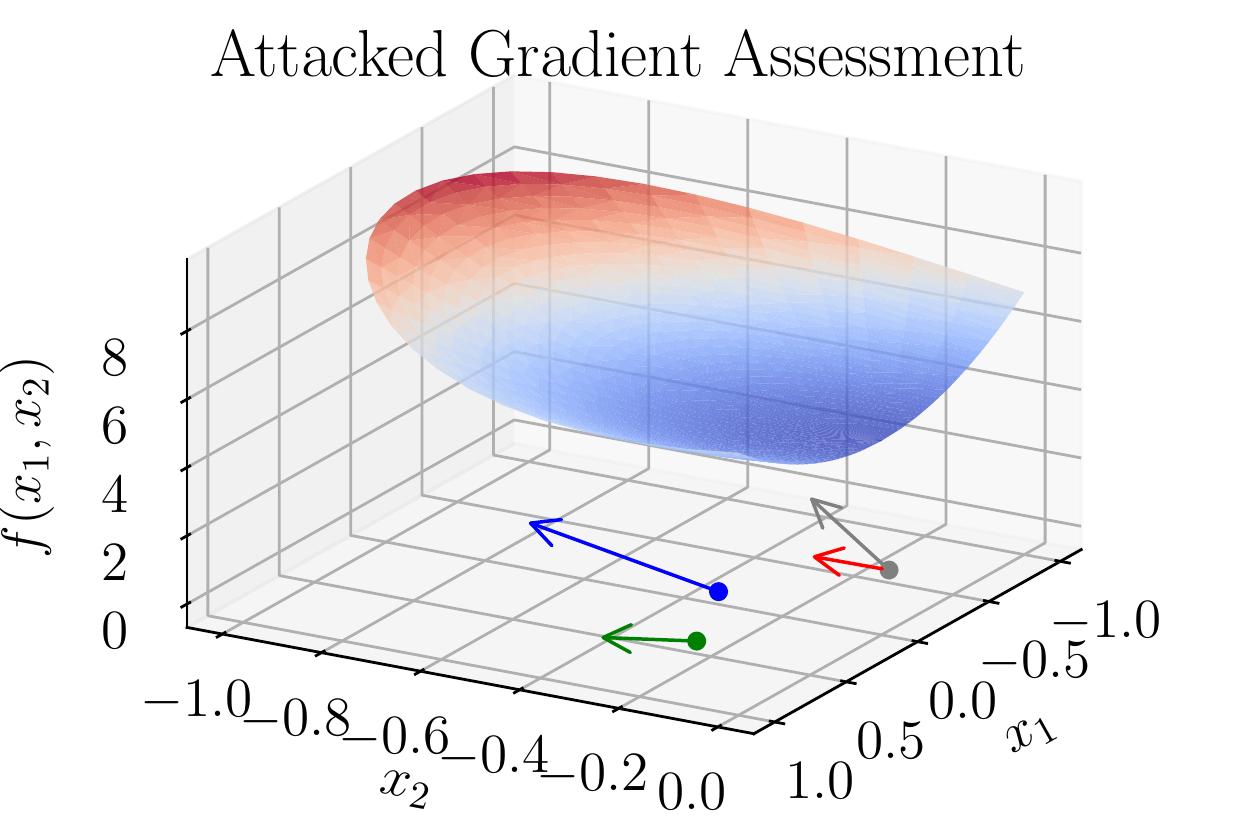}}
	\ifarticle
	
	\fi
	\subfloat[]{\includegraphics[scale = 0.5]{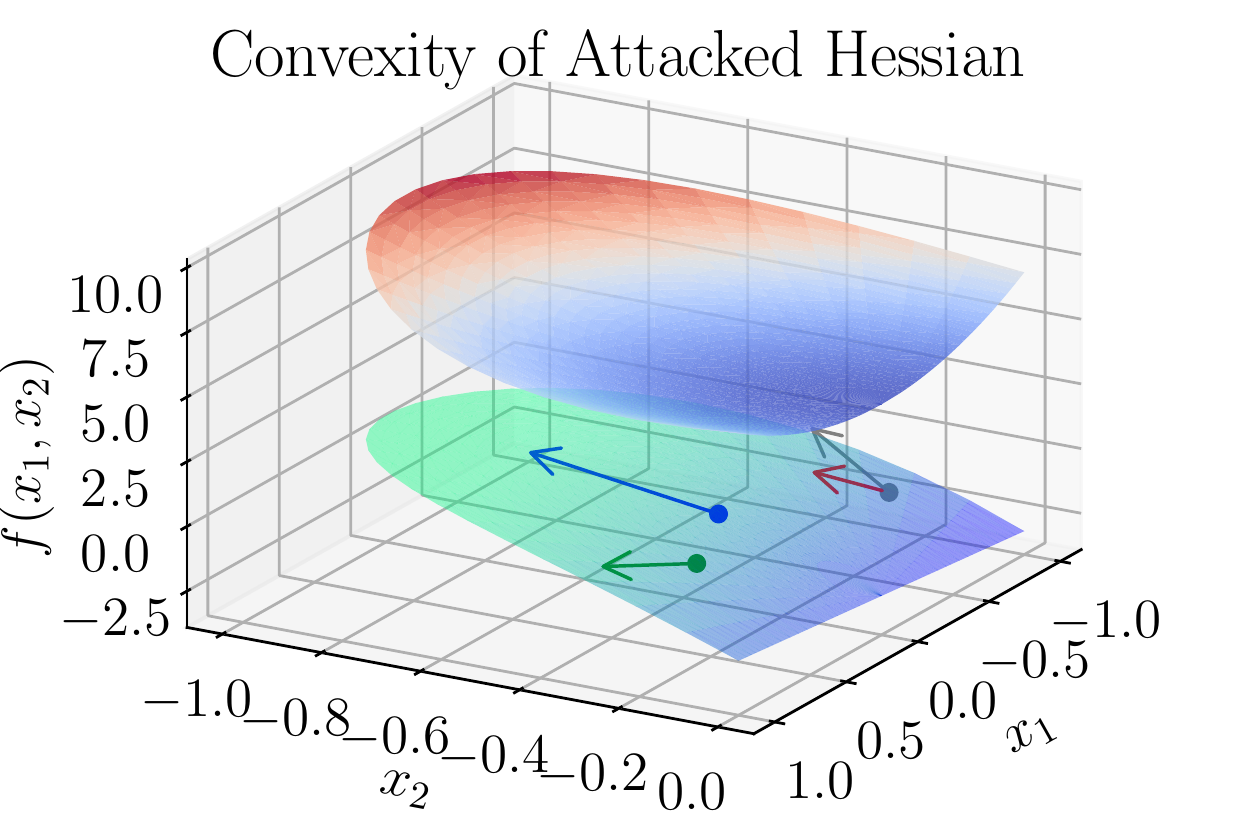}}
	\caption{The introduction of noise shifts the gradients computed using the above algorithm, shown here as shifting the unattacked (gray) vector to a new (attacked, red) position. It can be readily seen that a sufficient displacement will cause the estimated curvature (Hessian) to become nonconvex, as shown by the inferred surface represented in Subfigure (b).  To avoid detection, an attacker would thus need to use a very small injected signal, but this would limit their ability to stop convergence during the attack.}
	\label{fig:attackedVectors}
\end{figure}

\subsubsection{Assess Convexity}
For $f(x)$ to be convex, the Hessian must be positive semi-definite, i.e. $\lambda_{\text{min}}(H) \geq 0$ where  $\lambda_{\text{min}}(H)$ is the least eigenvalue of $H$.

After solving for a local approximation of the Hessian as above, the eigenvalues of the Hessian are evaluated. If the Hessian is not found to be positive semi-definite, then we suspect the $x$-actor may be injecting noise into his updates, as shown in Figure \ref{fig:attackedVectors}.

\subsubsection{Misclassifications}
\begin{figure}
	\centering
	\subfloat[]{\includegraphics[scale = 0.5]{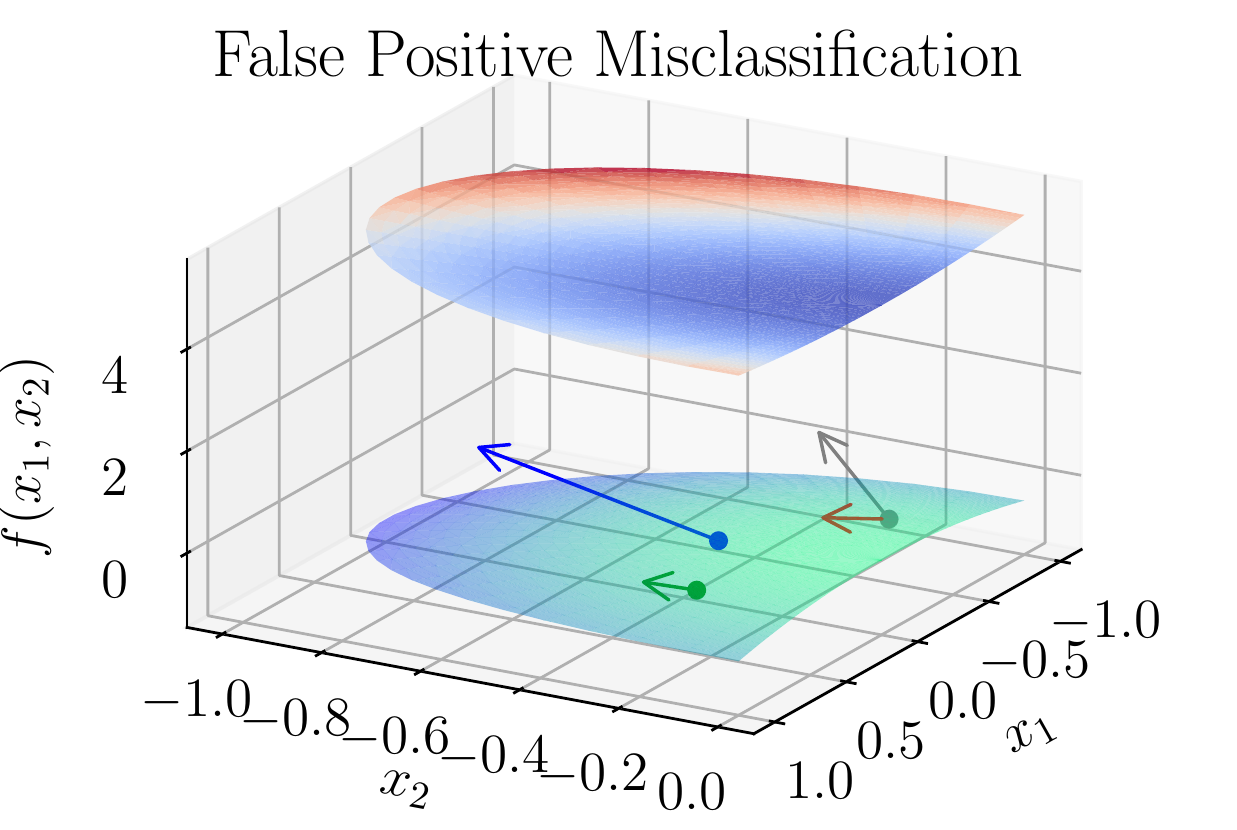}}
	\ifarticle
	
	\fi
	\subfloat[]{\includegraphics[scale = 0.5]{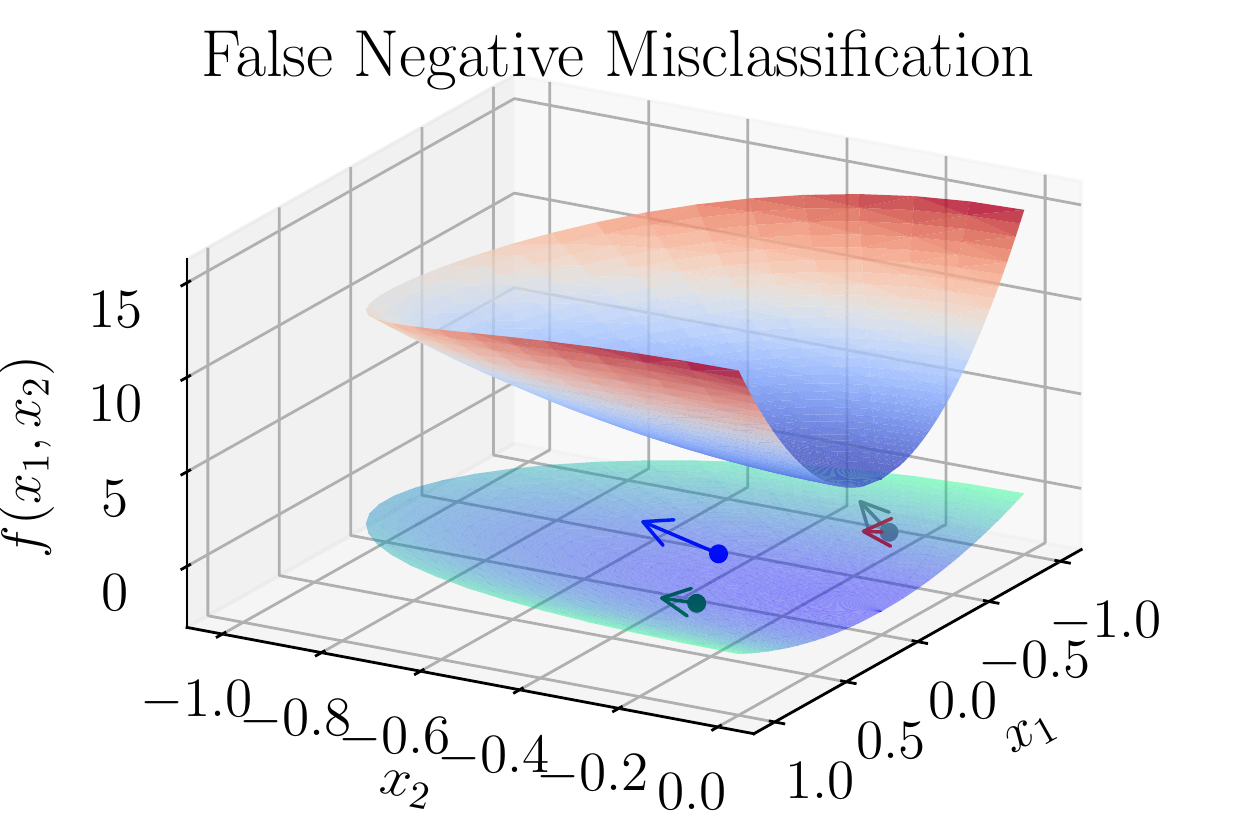}}
	\caption{Depiction of misclassification errors in the noise detection algorithm. In weakly convex regions of a function, numeric conditioning errors can lead to false positives (subfigure a), whereas problems which have a very strong gradient may converge before the noise detection algorithm is able to identify the injection of noise into the algorithm (subfigure b). In both cases, the scheme for selecting points used to assess the gradient can lead to more accurate assessment of the Hessian, reducing misclassification errors.}
	\label{fig:misclassification}
\end{figure}

As attack detection can be viewed as a classification problem, it may be subject to two possible errors, depicted in Figure \ref{fig:misclassification}.

\begin{itemize}
	\item False Positives (Type 1 errors): These occur when an attack is detected, even though no attack took place. This can occur when numeric issues create a slight local non-convexity during solution iterations, even though the underlying objective function is, in fact, convex. These numeric issues can result from the x-update, y-update, gradient calculation, or Hessian calculation, but are aggravated when the algorithm is close to the objective function and the system of equations used to generate the Hessian is poorly conditioned.
	\item False Negatives (Type 2 errors): In this type of error the algorithm fails to detect an attack, because the magnitude of the attack was not sufficient to create a non-convexity at the tested points. A small number of false positives are unavoidable: some attacks are not able to prevent convergence within a small number of iterations, and the attack detection algorithms may never be able to collect enough points to detect an attack. For problems that converge with a moderate number of iterations under attack, there may not be a combination of iterate points which result in a nonconvex Hessian. For other cases, lowering the false negative rate requires testing a variety of point combinations in order to maximize the likelihood of finding a point where the finite-difference Hessian estimate is nonconvex.
\end{itemize}

\section{Simulation and Results}

To verify the effectiveness of this problem under a variety of conditions, we simulate a large number of sample problems (both attacked and unattacked) and evaluate the effectiveness of the attack detection algorithm.

We have chosen to consider quadratic programs (QPs): they are widely used in power systems research (e.g. optimal power flow \cite{Kekatos2013}, optimal dispatch \cite{LeFloch2016a, Dobbe2016}), subsume the set of linear programming problems, and yet offer many computational benefits (efficient off-the-shelf solvers, easy analytic solutions, easily visualized). In the future, we would like to extend this to other classes of convex problems (e.g. SOCP, SDP), but simulating hundreds of those problems would be computationally burdensome.

\subsection{ Problem: Random Quadratic Programs}
We consider unconstrained quadratic programs of the form

\begin{align*}
\min_{x,z} &\quad x^T P x + c^T x + z^T Q z + d^T z\\
\text{s.t.} &\quad Ax + Bz = c
\end{align*}

For simplicity, we examine problems where $c=0$ and $A$ and $B$ have a single non-zero value per row, as described below. This can be thought of as consensus problems where some cost information is private.

\subsubsection{Problem Generation}
Problems are generated randomly as follows: 
\begin{enumerate}
\item The size of the problem is determined by randomly choosing the magnitudes $m,n,p$. The user defines a maximum size $\texttt{maxdim}$ and $n$ and $m$ are drawn as integers on the uniform distribution $[1, \texttt{maxdim}]$. The number of consensus constraints $p$ is then randomly chosen by drawing an integer on the interval $[1,\text{min}(m,n)]$. 
\item The problem will then be characterized by: 
\begin{itemize}
\item Variables $x\in\mathbb{R}^n,z\in\mathbb{R}^m$, 
\item Cost terms $P \in \mathbb{R}^{n \times n}, Q \in \mathbb{R}^{m \times m},c\in\mathbb{R}^n,d\in\mathbb{R}^m$, 
\item Constraint matrices $A\in\mathbb{R}^{p \times n}, B\in\mathbb{R}^{p \times m}$
\end{itemize}
\item Quadratic cost matrices $P$ and $Q$ are created by first randomly composing $L_P \in \mathbb{R}^{n \times n}, L_Q \in \mathbb{R}^{m \times m}$ with entries drawn randomly from  the uniform distribution on $[-S, S]$. $P$ and $Q$ are then constructed to be symmetric positive semi-definite by construction as $P = (L_P)^T L_P, Q = (L_Q)^T L_Q$. To confirm that these are positive semi-definite, the eigenvalues are computed and checked to be greater than 0.
\item $c$ and $d$ are generated by randomly drawing values from $[-S^2, S^2]$, resulting in entries that are the same magnitude as those in $P$ and $Q$. 
\item A is composed as 
$A = [0 \in \mathbb{R}^{(n-p)\times(n-p)}, I \in \mathbb{R}^{p \times p}]$ and B is composed as $B = [I \in \mathbb{R}^{p \times p}, 0 \in \mathbb{R}^{(m-p)\times(m-p)}]$
\end{enumerate}

\subsection{Solution Method}
We consider the case where consensus is desired between two actors who do not want to share information about their cost information $P,c$ and $Q,d$. In this scenario, ADMM is a well-suited tool for achieving optimality without compromising privacy.

For this problem, the ADMM algorithm can be stated as:
\begin{align*}
x^{k+1} &= \text{argmin}_x \quad x^T P x + c^T x + \frac{\rho}{2}||Ax+Bz^{k}- c + u^k||^2_2 \\
z^{k+1} &= \text{argmin}_z \quad z^T Q z + d^T z + \frac{\rho}{2}||Bz + Ax^{k+1} - c + u^k||^2_2 \\
u^{k+1} &= u^k + Ax^{k+1} + Bx^{k+1} - c
\end{align*}

The x-update and z-update steps constitute unconstrained QPs, and can be solved analytically as described in the Appendix. These analytic update steps were used to avoid rounding issues resulting from using an iterative numeric solver.

\subsection{Attack Implementation}
Attack was simulated by multiplying the optimal $x$-update values by a vector with entries randomly selected from ${+10\%, -10\%}$, i.e. a Bernoulli distribution, at each iteration. It was found that choosing uniformly distributed zero-centered noise was not sufficient to prevent convergence, as small-magnitude noise entires allow the algorithm to converge faster than the attacker is able to delay convergence.

A set of 10,000 QPs was simulated without attack, and were found to converge in a median of 63 iterations, with a mode of 16 iterations and a maximum value of 216 iterations. When simulated with the attack described above,  convergence was significantly hindered on the same QPs: 86\% were not able to converge in 300 iterations, and 87\% do not converge in 500 iterations. The 14\% of QPs which converge in less than 300 iterations were found to behave similarly to the unattacked QPs, likely due to a set of costs which are strongly convergent even in the presence of attack. 

\begin{figure}
	\centering
	\ifarticle \newcommand{\figwidth}{1.0} \else \newcommand{\figwidth}{0.75} \fi
	\includegraphics[width=\figwidth\columnwidth]{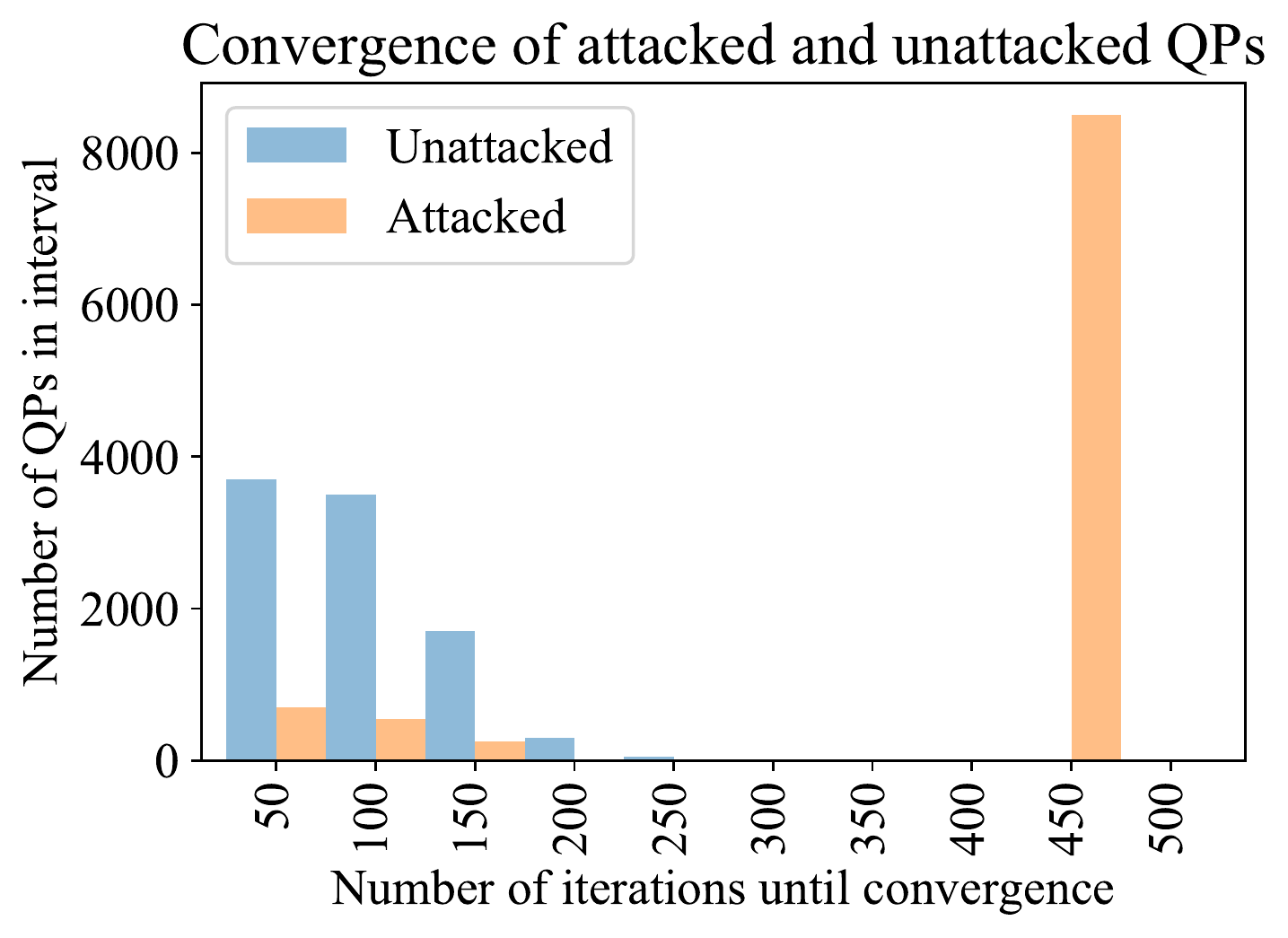}
	\caption{Convergence results for 10,000 QPs, simulated both with and without noise-injection attack. Without attack, all simulations converge in less than 300 iterations. Under attack, 86\% of problems do not converge in 500 iterations (at which time calculation was stopped).}
	\label{fig:qpconvergence}
\end{figure}

The distribution of the iterations until convergence for both the unattacked and attacked QPs is depicted in Figure \ref{fig:qpconvergence}.

\subsection{Attack Detection Implementation}
We implement the algorithm described above. As we initially consider $p\leq2$, we need to complete at least 4 iterations in order to create two sets of difference equations.

To improve conditioning, for iteration $i$ we consider the difference equations with $i$ and the reference point, and compare with iteration 2 and iteration $\lfloor i / 2 \rfloor$. Using the $i, i-1, i-2$ iterates was found to result in poor conditioning and a very high (17\%) false positive rate, though a slightly lower false negative rate.

To ensure that we can solve for the values of the Hessian entries, we check that the difference vectors are not collinear.

\subsection{Results}

We show results for 1000 simulated problems, of which 500 were attacked and 500 were unattacked. The results are presented as a \textit{confusion matrix} in Table \ref{tab-results}, which shows both accurate classifications and misclassifications.  The number of problems of each type is shown in Table \ref{tab-results}\subref{tab-results-abs}; the same results are shown in proportional form in Table \ref{tab-results}\subref{tab-results-rel}.

\begin{table}[]
	\centering
	\caption{Results for a simulation of attack detection. 500 quadratic programs were generated without attack, and false positives (upper right quadrant) were identified. An additional 500 quadratic programs were generated with simulated attack, and false negatives measured.}
	\label{tab-results}
	\subfloat[Results, total number of simulations]{
		\label{tab-results-abs}
	\begin{tabular}{|l|c|c|}
		\hline
		& Attacked & Unattacked  \\ \hline
		Attack Detected 		&   479     &  0         \\ \hline
		No Attack Detected 	 &   21     &  500     \\ \hline \hline
		\textit{Total} 				& \textit{500}   & \textit{500} \\ \hline
	\end{tabular}}

	\subfloat[Results, simulation error rate]{
		\label{tab-results-rel}
	\begin{tabular}{|l|c|c|}
	\hline
	& Attacked & Unattacked  \\ \hline
	Attack Detected 		&   0.958     &  0         \\ \hline
	No Attack Detected 	 &   0.042     &  1.0     \\ \hline
	\textit{Total} 				& \textit{1.0}   & \textit{1.0} \\ \hline
\end{tabular}}

\end{table}

We consider two types of errors: Type 1 (false positive, upper right quadrant) errors, and Type 2 (false negative, bottom leftquadrant ) errors. The experimental results were compared when computing the local update with both numeric solvers and with an analytic solution; both approaches were found to give comparable results. 

The choice of iterates used to compute the difference equations was found to have a significant impact on the results: as only $n+1$ points are needed but $k$ iterates are available, there are $\frac{k!}{(n+1)!(k-n-1)!}$ possible choices for points to use, assuming that the most recent point is chosen as a base. A naive approach in which the $n+1$ most recent iterates were used was found to result in poor numeric conditioning and a large number of Type 1 errors; the results presented here select the $n+1$ iterates to be evenly spaced between iterate 2 and iterate $k$.

\section{Limitations}
This algorithm has been tested with constrained and unconstrained QPs of relatively small dimension, and in development was also tested with linear programs.  We have not explored other problem types (e.g. second-order cone programs, semi-definite programs) where the Hessian may change over the iterates. We also have not expressly addressed the theoretical implications of using this on strongly-convex problems (e.g. problems which can be bounded below by a QP) compared with weakly-convex problems.

\subsection{High-Dimensional Problems}
We have previously assumed that $n+1$ iterates are needed to assess the validity of the $x$-update step. Where $n$ is large, collecting sufficient iterates may be costly in time or computation resources. Further, the problem may converge before $n+1$ iterates are collected (Note, this is not a problem if the goal is to prevent convergence-stopping attacks).

We consider two scenarios: one where the full set of x-update variables are made public, and one where only the variables associated with the linking constraint are made public.

\subsubsection*{All x-variables are public}
When all x-variables are made public, the gradient is constructed as:
\begin{align*}
A & \in \mathbb{R}^{p\times n} \\
x_{\text{pub}} & \in \mathbb{R}^n \\
(z^i-z^{i-1}) &\in \mathbb{R}^m \\
(y^i - \rho B(z^i-z^{i-1}) &\in \mathbb{R}^p \\
-A^T (y^i - \rho B(z^i-z^{i-1})) &\in \mathbb{R}^n 
\end{align*}

The resulting gradient contains all $x$-dimensions, but will be rank deficient, having span $\mathbb{R}^p$ and nullspace $n-p$. In this case, each new point gives us $n$ new dimensions, but we also seek to determine the Hessian in $n$ dimensions and $n^2$ unknowns. We thus need to gather $n+1$ points.

\subsubsection*{Only linking dimensions are public}
If only linking dimensions are public, the gradient is constructed as: 
\begin{align*}
A & \in \mathbb{R}^{p\times p} \\
x_{\text{pub}} & \in \mathbb{R}^p \\
-A^T (y^i - \rho B(z^i-z^{i-1})) &\in \mathbb{R}^p 
\end{align*}
In this case, the gradient estimate fully spans the linking dimensions, and each new point gives us $p$ new equations. We seek the Hessian in $p$ dimensions and $p^2$ unknowns, thus needing to gather $p+1$ datapoints. 

Either way, we only are able to assess nonconvexity in the linking dimensions, and are unsure of activity in the other dimensions.

\section{Extensions}

\subsection{Reducing Error Rates}
As described above, the false positive and false negative rates can be reduced by tuning the algorithm:
\begin{itemize}
	\item \textit{Reducing False Positives:} As previously described, false positives result from poor numeric conditioning in computing the Hessian from the system of difference equations. Reducing these errors requires testing the condition number of the system of equations, and potentially also testing the condition number of the resulting Hessian estimate. Neither topic has been explored here. 
	\item \textit{Reducing False Negatives:} False negatives occur when the tested set of points do not indicate a non-convexity; this is most likely to occur when the magnitude of attack is small relative to the curvature of the $x$-update function. Reducing this error rate requires testing a variety of combinations of points in order to maximize the likelihood of finding a point where magnitude of the attack is large relative to the local curvature of the function.
\end{itemize}

\subsection{Improving Algorithmic Efficiency}
We have not explored opportunities to improve the efficiency of the algorithm, but highlight three opportunities here:
\begin{itemize}
	\item \textbf{Exploiting symmetry in Hessian:} We have treated the Hessian as having $n^2$ unknowns, but as it is symmetric there are actually only $n(n+1)/2$ independent terms. Exploiting this structure would mean that only $\lceil (n+1)/2\rceil +1 $ points are needed, but also requires pruning the equation set to avoid over-determining the system of equations.
	\item \textbf{Exploiting problem-specific structure of Hessian:} For quadratic programs, the Hessian matrix is constant throughout the problem, meaning that the difference equations do not need to be constructed with respect to a single point, but rather can be constructed using all combinations of points. In this case, we need to only use $j$ iterates such that ${j \choose 2} > \lceil(n+1)/2\rceil + 1$, allowing the algorithm to be started faster, and allowing for more robust checking in the case of poor conditioning.
	\item \textbf{Choosing point set:} If the Hessian is computed multiple times per iteration to reduce error rates as described above, choosing the iterates wisely can improve numeric conditioning and reduce the need for extra computations. This has not been explored.
\end{itemize}

\section{Extension: Fully-decentralized optimization}

As the central information hub, the aggregator in a decentralized optimization problem must be trusted to provide fair updates to all the nodes. When the aggregator and the local nodes have the same incentives -- e.g., they are all computing resources owned by the same entity -- this is unlikely to present a conflict of interests. 

However, when the aggregator has different incentives than the compute nodes, the local nodes may not trust the central aggregator, e.g. if the aggregator can extract additional profits by acting as a market maker.  Further, aggregator-based decentralized computation is subject to several other weaknesses: it requires a high communication overhead by coordinating message-passing between all nodes (a significant hurdle for highly distributed systems like energy devices), and introduces a central point of failure in the case of communication/power outage or cyberattack.

To address these issues with aggregator-coordinated decentralized optimization, \textit{fully-decentralized} algorithms have been developed which take advantage of problem structure to achieve consensus between nodes that directly share constraints, rather than passing information from all nodes to a centralized aggregator \cite{Mota2012}. 

\subsection{Fully-Decentralized ADMM}
As an example, the ADMM algorithm from above can be expressed in fully-decentralized form by introducing local copies $u_x$ and $u_z$ of the penalty variables:
\begin{align}
x^{k+1} &:= \; \text{argmin}_{x \in \mathcal{X}} \; f(x) + \frac{\rho}{2}\lVert Ax + Bz^k - c + u_x^k \rVert^2_2 \\
z^{k+1} &:= \; \text{argmin}_{z \in \mathcal{Z}} \; g(z) + \frac{\rho}{2}\lVert Ax^{k+1} + Bz - c + u_z^k \rVert^2_2 \\
u_x^{k+1} &:= \; u_x^k + Ax^{k+1} + Bz^{k+1} - c\\
u_z^{k+1} &:= \; u_z^k + Ax^{k+1} + Bz^{k+1} - c
\end{align}

\ifarticle \newcommand{\figwidth}{0.75} \else \newcommand{\figwidth}{1.0} \fi

\begin{figure}
	\begin{tikzpicture}[
	scale=\figwidth,
	every node/.style={transform shape},
	bluenode/.style={circle, draw=blue!60, fill=blue!5, very thick, minimum size=7mm},
	rednode/.style={circle, draw=red!60, fill=red!5, very thick, minimum size=5mm},
	]
	
	\node[bluenode, align=center]   (leftend)[text width=1.7cm]   {upstream nodes};
	\node[bluenode]  (x)       [right=of leftend] {$ f(x), x\in\mathcal{X}$};
	\node[bluenode]   (z)       [right=of x] {$g(z), z\in\mathcal{Z}$};
	\node[bluenode, align=center]  (rightend)       [right=of z,text width=1.7cm] {downstream nodes};
	
	\draw [->]  ([yshift= 0.2cm]  leftend.east) -- node[yshift=0.3cm]   {...}        ([yshift=  0.2cm]  x.west);
	\draw [<-]  ([yshift=-0.2cm] leftend.east) -- node[yshift=-0.3cm]  {$x^k$} ([yshift=-0.2cm]  x.west);
		
	\draw [->]  ([yshift= 0.2cm]  x.east) -- node[yshift=0.3cm]   {$x^k$}        ([yshift=  0.2cm]  z.west);
	\draw [<-]  ([yshift=-0.2cm] x.east) -- node[yshift=-0.3cm]  {$z^k$} ([yshift=-0.2cm]  z.west);
	
	\draw [->]  ([yshift= 0.2cm]  z.east) -- node[yshift=0.3cm]   {$z^k$}        ([yshift=  0.2cm]  rightend.west);
	\draw [<-]  ([yshift=-0.2cm] z.east) -- node[yshift=-0.3cm]  {$\ldots$} ([yshift=-0.2cm]  rightend.west);
	\end{tikzpicture}
	\caption{Fully-decentralized optimization problem structure, highlighting two nodes and noting how this can be indefinitely expanded with the addition of further upstream and downstream nodes. Each node holds private information on its own objective function, constraints, and dual variable (e.g. $u_x^k$), and accepts updates from its immediate neighbors.}
	\label{fig:nodefigure_fullyDecentOpt}
\end{figure}
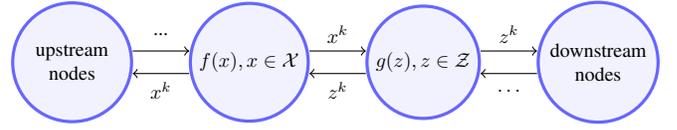

Under restrictions described in \cite{Mota2012}, this can be shown to have the same convergence and optimality guarantees as conventional aggregator-based systems.

In systems where nodes are very sparsely connected, this can produce a significant reduction in  communication overhead- e.g. in power systems we may seek an optimum amongst millions of nodes, but each node is only connected to its parent and one or two downstream nodes. Rather than requiring that an aggregator handle millions of connections, nodes pass messages with their neighbors, ultimately bringing the full system into optimum. For examples of this approach applied to power systems, see e.g. \cite{Sulc2014,Tsai2015,Wang2016}.

\section{Attack Vectors in Fully-Decentralized Optimization}

As nodes in a fully-decentralized network are only connected with their neighbors, they must assess whether updates represent the true state of the network, or whether the neighbor has been compromised and the update is spurious. The following sections highlight the unique challenges of detecting, localizing, and mitigating attacks in a fully-decentralized system, as shown graphically in Figure \ref{fig:nodeAttackPossibilities}.

\ifarticle \renewcommand{\figwidth}{0.70} \else \renewcommand{\figwidth}{1.0} \fi

\begin{figure}
	\centering
	\subfloat[]{
		\begin{tikzpicture}[
		scale=\figwidth, every node/.style={transform shape},
		bluenode/.style={circle, draw=blue!60, fill=blue!5, very thick, minimum size=7mm}, rednode/.style={circle, draw=red!60, fill=red!5, very thick, minimum size=5mm}, ]
		\node[bluenode] (w)      [right=of leftend] {$\begin{gathered}w^k = \text{argmin}_w [ \\L_{\rho}(h(w),x^k, u^k_w) ]\end{gathered}$}; 
		\node[bluenode]  (x)       [right=of w] {$\begin{gathered}x^k = \text{argmin}_x [ \\L_{\rho}(f(x),w^k,z^k, u^k_x) ]\end{gathered}$}; 
		\node[bluenode]   (z)       [right=of x] {$\begin{gathered}z^k = \text{argmin}_z [ \\L_{\rho}(g(z),x^{?k}, u^k_z) ]\end{gathered}$}; 
		
		\draw [->]  ([yshift= 0.2cm]  w.east) -- node[yshift=0.3cm]   {$w^k$}        ([yshift=  0.2cm]  x.west);
		\draw [<-]  ([yshift=-0.2cm] w.east) -- node[yshift=-0.3cm]  {$x^k$} ([yshift=-0.2cm]  x.west);
		
		\draw [->]  ([yshift= 0.2cm]  x.east) -- node[yshift=0.3cm]   {$x^{?k}$}        ([yshift=  0.2cm]  z.west);
		\draw [<-]  ([yshift=-0.2cm] x.east) -- node[yshift=-0.3cm]  {$z^k$} ([yshift=-0.2cm]  z.west);
		\end{tikzpicture}
		\label{fig:attackscenarios_noattack}
	}
	
	\subfloat[]{
		\begin{tikzpicture}[
		scale=\figwidth, every node/.style={transform shape},
		bluenode/.style={circle, draw=blue!60, fill=blue!5, very thick, minimum size=7mm}, rednode/.style={circle, draw=red!60, fill=red!5, very thick, minimum size=5mm}, ]
		\node[bluenode] (w)      [right=of leftend] {$\begin{gathered}w^k = \text{argmin}_w [ \\L_{\rho}(h(w),\tilde{x}^k, u^k_w) ]\end{gathered}$};
		\node[rednode]  (x)       [right=of w] {$\begin{gathered}\tilde{x}^k \not= \text{argmin}_x [ \\L_{\rho}(f(x),w^k,z^k, u^k_x) ]\end{gathered}$};
		\node[bluenode]   (z)       [right=of x] {$\begin{gathered}z^k = \text{argmin}_z [ \\L_{\rho}(g(z),x^{?k}, u^k_z) ]\end{gathered}$};
		
		\draw [->]  ([yshift= 0.2cm]  w.east) -- node[yshift=0.3cm]   {$w^k$}        ([yshift=  0.2cm]  x.west);
		\draw [<-]  ([yshift=-0.2cm] w.east) -- node[yshift=-0.3cm]  {$\tilde{x}^k$} ([yshift=-0.2cm]  x.west);
		
		\draw [->]  ([yshift= 0.2cm]  x.east) -- node[yshift=0.3cm]   {$x^{?k}$}        ([yshift=  0.2cm]  z.west);
		\draw [<-]  ([yshift=-0.2cm] x.east) -- node[yshift=-0.3cm]  {$z^k$} ([yshift=-0.2cm]  z.west);
		\end{tikzpicture}
		\label{fig:attackscenarios_xattack}
	}
	
	\subfloat[]{
		\begin{tikzpicture}[
		scale=\figwidth, every node/.style={transform shape},
		bluenode/.style={circle, draw=blue!60, fill=blue!5, very thick, minimum size=7mm}, rednode/.style={circle, draw=red!60, fill=red!5, very thick, minimum size=5mm}, ]
		\node[rednode] (w)      [right=of leftend] {$\begin{gathered}\tilde{w}^k \not= \text{argmin}_w [ \\L_{\rho}(h(w),x^k, u^k_w) ]\end{gathered}$};
		\node[bluenode]  (x)       [right=of w] {$\begin{gathered}x^k = \text{argmin}_x [ \\L_{\rho}(f(x),\tilde{w}^k,z^k, u^k_x) ]\end{gathered}$};
		\node[bluenode]   (z)       [right=of x] {$\begin{gathered}z^k = \text{argmin}_z [ \\L_{\rho}(g(z),x^{?k}, u^k_z) ]\end{gathered}$};
		
		\draw [->]  ([yshift= 0.2cm]  w.east) -- node[yshift=0.3cm]   {$\tilde{w}^k$}        ([yshift=  0.2cm]  x.west);
		\draw [<-]  ([yshift=-0.2cm] w.east) -- node[yshift=-0.3cm]  {$x^k$} ([yshift=-0.2cm]  x.west);
		
		\draw [->]  ([yshift= 0.2cm]  x.east) -- node[yshift=0.3cm]   {$x^{?k}$}        ([yshift=  0.2cm]  z.west);
		\draw [<-]  ([yshift=-0.2cm] x.east) -- node[yshift=-0.3cm]  {$z^k$} ([yshift=-0.2cm]  z.west);
		\end{tikzpicture}	
		\label{fig:attackscenarios_wattack}
	}
	
	\caption{Potential attack scenarios in a fully-decentralized optimization scenario, where the $z$-update node is attempting to establish the veracity of a received update $x^{?k}$. Without knowing details of the private objective functions $f(x),h(w)$ and constraints $x\in\mathcal{X}, w\in\mathcal{W}$ it is difficult to detect and localize (or mitigate) an attack. The $z$-update node is not generally able to determine the difference between the unattacked scenario shown in \protect\subref{fig:attackscenarios_noattack}, an attack by the immediate upstream neighbor $\tilde{x}^k$ as in \protect\subref{fig:attackscenarios_xattack}, and an attack by a node further upstream such as $\tilde{w}^k$ as shown in \protect\subref{fig:attackscenarios_wattack}. }
	\label{fig:nodeAttackPossibilities}
\end{figure}
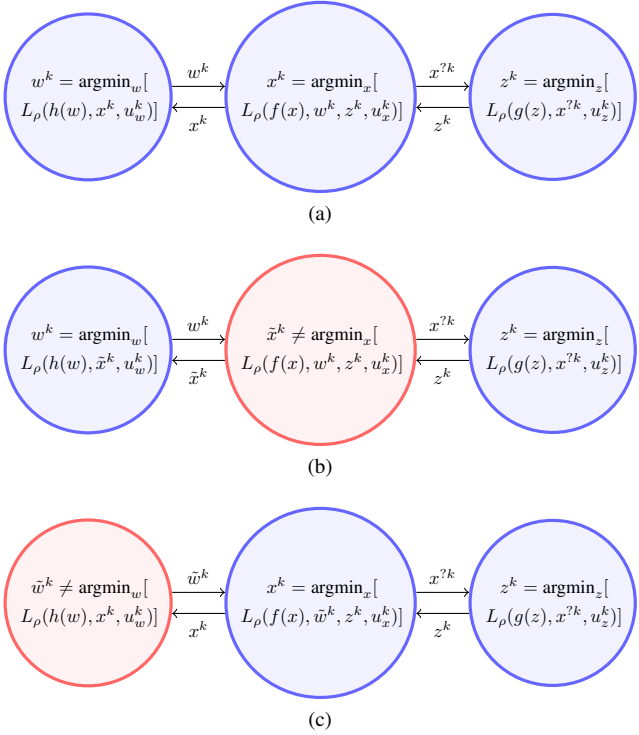

\subsection{Attack Vector: Private Infeasibility Attack}
In Section \ref{sec:attackvector_overview}, we highlighted how a malicious node may distort its private constraints to shift the equilibrium out of the operationally feasible region.

This attack is not identifiable in a fully-decentralized system, as a node can not in general discern between the update corrupted by its immediate neighbor

$$x^{? k} \overset{?}{=} \text{argmin}_{x\in\tilde{\mathcal{X}}} L_{\rho}(x, z^{k-1}, w^{k-1},u_x^{k-1}) $$ 
and a best response from a neighbor in reaction to a corrupted signal from the upstream $h(w)$ node:
$$x^{? k} \overset{?}{=} \text{argmin}_{x\in\mathcal{X}} L_{\rho}(x, z^{k-1},\tilde{w}^{k-1}, u_x^{k-1})$$

Alternately, on seeing $x^{?k} \in \{\mathcal{X}_{\text{pub}} | z^{k-1}\}$ it is necessary to know $w^k$ in order to assess whether the problem is truly feasible given the \textit{full} state of the system, $x^{?k} \overset{?}{\in} \{\mathcal{X}|w^{k-1}, z^{k-1} \}$ (this can be extended to more complex system architectures with more upstream/downstream nodes).

However, in a fully-decentralized system, $w$ is not generally visible to $x$ and so it is not possible for the $z$-update node to assess whether $x^{?k} \in \{\mathcal{X}_{\text{pub}}|w^k\}$.

This means that detection, localization, and mitigation are all not possible in a conventional fully-decentralized system, as a global information layer is needed to assess the feasibility of the received updates, identify nodes which may be causing infeasibility, and construct a best response.

\subsection{Attack vector: Infeasible Linking Constraint}\label{sec:attack1}

Similarly, on receiving an update $x^{?k}$ which would violate the linking constraint $\{Ax^{?k} + Bz = c | z \in \mathcal{Z}\}$ the $z$-update node is not able to discern between an infeasible update created by a neighboring node:
$$x^{? k} \overset{?}{=} x^{\star k} + \varepsilon^k$$
and an infeasible signal resulting from a malicious upstream/downstream node: 
 
 \begin{align*}
x^{? k} &\overset{?}{=} \text{argmin}_{x\in\mathcal{X}} L_{\rho}(x, z^{k-1},\tilde{w}^{k-1}, u_x^{k-1})\\
&= \text{argmin}_{x\in\mathcal{X}} L_{\rho}(z^{k-1},z,w^{\star k-1} + \varepsilon, u_x^{k-1})
 \end{align*}

However, if information on public constraints $\mathcal{Z}_\text{pub}, \mathcal{X}_\text{pub}$ is knowable to all participants, these constraints can be integrated into the local optimization problems as additional (publicly known) constraints, converting the update step from 
$$
x^{k+1} = \text{argmin}_{x\in\mathcal{X}} L_{\rho}(x,z^k, u^k)
$$
to 
$$
x^{k+1} = \text{argmin}_{x\in \mathcal{X} \cap \{\mathcal{Z}_{\text{pub}}|A^T x+B^T z^k-c\}} L_{\rho}(x,z^k, u^k)
$$
This would mean that \textit{any} update which violates this constraint must be the result of a malicious node, and the problem can be localized.

\subsubsection*{Attack Mitigation}
If a feasible solution exists, it must satisfy the constraints:
\begin{align*}
x^{\star} \in \mathcal{Z} \subset \mathcal{X}_{\text{pub}}\\
z^{\star} \in \mathcal{Z} \subset \mathcal{Z}_{\text{pub}}\\
A^T x^{\star} + B^T z^{\star} = c
\end{align*}
Therefore, at each step we can project any iterate onto the feasible set described by the constraints above:
$$
\hat{x}^k = \text{argmin}_{z\in \{\mathcal{Z}|A^T x + B^T z = c \} } ||x^{? k} - z||_2^2
$$
In an unattacked case this projection simply accelerates convergence of the standard ADMM algorithm.

In the case of attack, projection creates a $\hat{x}^k \in \{\mathcal{Z}|A^T x + B^T z = c \}$ which is the \textit{best response} to the attack $\tilde{x}^k$. While this best response will in general be suboptimal, it is feasible and will let the primal residual $r = A x^k + B z^k - c$ converge to zero even in the presence of attack.

\subsection{Attack Vector: Zero-Mean Noise Injection}

Because it is only dependent on local information (its own history and updates received from the neighbors), each node in a fully-decentralized network is able to carry out the detection strategy outlined in Section \ref{sec:algorithm_noisedetection}  for identifying the presence of a noise-injection attack.

However, in a fully-decentralized network the node would be challenged to identify whether the noise injection attack originates from an immediate neighbor, or from an upstream node. 

Specifically, if $h(w)$ represents the problem solved by the upstream node, the $z$-update node cannot differentiate between $x^{? k} \overset{?}{=} x^{\star k} + \delta(k)$ and 
\begin{align*}
x^{? k} &\overset{?}{=} \text{argmin}_{x\in\mathcal{X}} L_{\rho}(x, z^{k-1},\tilde{w}^{k-1}, u^{k-1})\\
&= \text{argmin}_{x\in\mathcal{X}} L_{\rho}(x,z^{k-1},w^{\star k-1} + \delta(k-1), u^{k-1})
\end{align*}
Nevertheless, although the malicious node cannot be localized, a zero-mean noise injection attack can be detected.

\section{Summary of Security Challenges}

Table \ref{tbl-security-feasibility}\subref{tbl-security-aggregator-coord} highlights the feasibility of achieving these goals in conventional aggregator-coordinated optimization, and Table \ref{tbl-security-feasibility}\subref{tbl-security-fully-decent} shows the same capabilities in a fully-decentralized system. We close by exploring approaches for reaching these goals in a fully-decentralized environment.

\begin{table}[h]
	\centering
	\caption{Security issues in aggregator-coordinated and fully-decentralized systems}
	\label{tbl-security-feasibility}
	\subfloat[Detection Feasibility on Aggregator-Coordinated System]{
		\begin{tabular}{|l|c|c|c|}
			\hline
			& Detect & Localize & Mitigate \\ \hline
			Private Infeasibility &   X     &  X        &   X       \\ \hline
			Linking Infeasibility &   X     &  X        &          \\ \hline
			Noise Injection       &    X    &   X       &          \\ \hline
		\end{tabular}
		\label{tbl-security-aggregator-coord}
	}
	
	\subfloat[Detection Feasibility on Fully-Decentralized System]{
		\begin{tabular}{|l|c|c|c|}
			\hline
			& Detect & Localize & Mitigate \\ \hline
			Private Infeasibility &         &       &          \\ \hline
			Linking Infeasibility &   X   &  X   &          \\ \hline
			Noise Injection       &    X   &       &          \\ \hline
		\end{tabular}
		\label{tbl-security-fully-decent}
	}
\end{table}

\subsection{Architectures for Security}

Although aggregator-based systems present some security and trust issues, the aggregator is able to take advantage of global information to check the validity of each node's updates, enabling detection, localization, and mitigation strategies, which are not available in the fully-decentralized system.

However, if the fully-decentralized system is augmented with a global information layer, the same security checks become possible in a fully-decentralized environment. In this section, we briefly discuss information architectures which could allow a fully-decentralized system to take advantage of the security benefits of an aggregator-coordinated system.

A number of different architectures might provide this security:
\begin{enumerate}
\item A centralized database held by a trusted authority, with security checks computed by each node \label{security:central}
\item A fully-connected network with securely signed messages, in which each node maintains a database of message history \label{security:fullconn}
\item A partially-connected network with \textit{bypass connections} to allow nodes to bypass suspect neighbors \label{security:bypass}
\item A decentralized peer-to-peer database with message histories, with security checks computed by each node \label{security:decentraldb}
\item A \textit{blockchain} used as a decentralized database to store the message histories, with \textit{smart contracts} used to compute security checks in a decentralized manner. \label{security:blockchain}
\end{enumerate}

\begin{figure}
	\centering
	\ifarticle \renewcommand{\figwidth}{1.0} \else \renewcommand{\figwidth}{0.75} \fi
	\includegraphics[width=\figwidth\columnwidth]{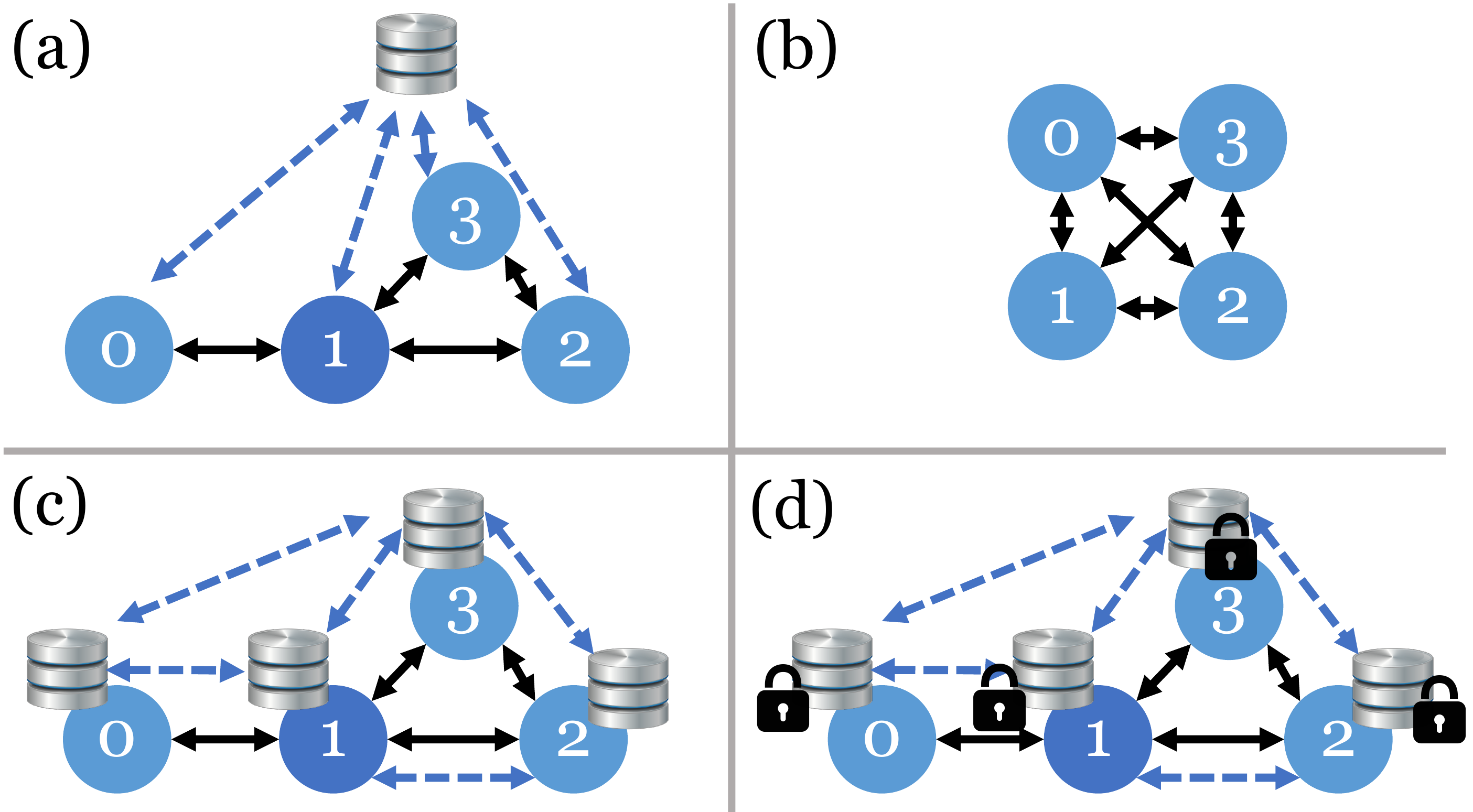}
	\caption{Potential architectures for allowing security checks for fully-decentralized optimization algorithms.}
	\label{fig:architectures}
\end{figure}

Option \ref{security:central} presents trust and monopoly distortion issues described above, in addition to a high communication overhead. Option \ref{security:fullconn} requires high computational overhead and will fail if there are dropouts in communication with part of the network; Option \ref{security:bypass} has gained some attention as a way to approach this while not requiring the same degree of communication redundancy. The decentralized database in Option \ref{security:decentraldb} is appealing, but requires a method for reconciling different versions of the database which might be proposed  by different neighbors. Further, both Option \ref{security:central} and Option \ref{security:decentraldb} do not guarantee that each node conducts the same security checks.

By contrast, \textit{blockchains} are decentralized databases which rely on securely signed messages and a consensus mechanism that simultaneously guarantees consistency across the network and provides secure timestamping.  \textit{Smart contracts} build on this architecture by guaranteeing the execution of simple computational functions as part of the consensus mechanism, thus guaranteeing consistent execution of security checks.

\begin{figure}
	\centering
	\ifarticle \renewcommand{\figwidth}{1.0} \else \renewcommand{\figwidth}{0.75} \fi
	\includegraphics[width=\figwidth\columnwidth]{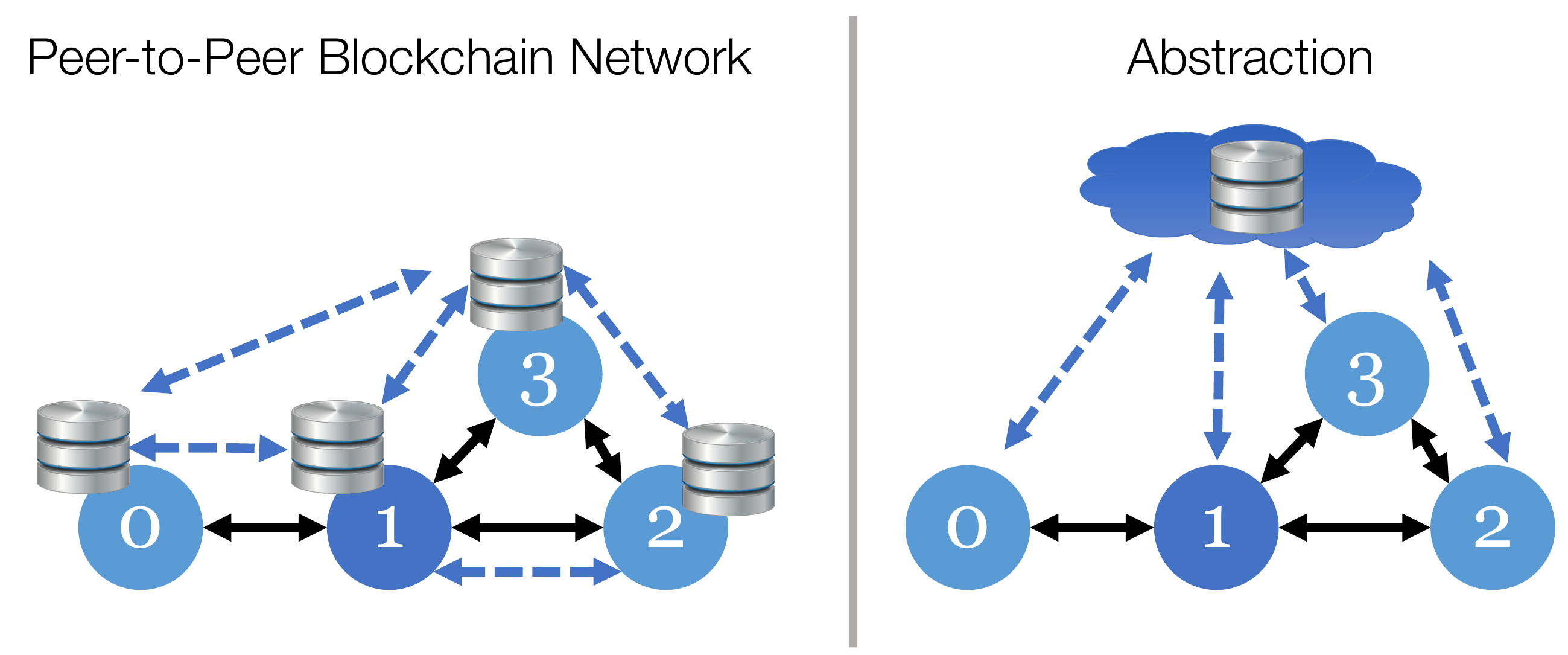}
	\caption{Illustration of how a blockchain-based network provides comparable security to an aggregator-coordinated architecture}
	\label{fig:blockchainAbstraction}
\end{figure}

Figure \ref{fig:blockchainAbstraction} conceptually shows how a blockchain-based system enables the same type of security checks as are offered by an aggregator system, by providing a cryptographically secured global information layer with guaranteed execution of the security checks outlined above.

Consequently, the features of cyberattack diagnostics algorithms for aggregator-coordinator systems can be accomplished in a fully decentralized setting, if a shared and secure information layer is provided, such as blockchains.


\section{Conclusion}
We have demonstrated the weaknesses of conventional distributed and fully-decentralized architectures to a number of attack vectors, including distortion of the private optimization problem, injection of noise, and distortion of the linking constraints in fully-decentralized networks. As optimization problems become distributed to consumer hardware (e.g. smartphones, smartmeters, autonomous vehicles, IoT devices), optimization algorithms becomes exposed to these threats along an exponentially larger attack surface. We examine methods for detecting and mitigating attacks. In particular, we detail a detection algorithm for noise-injection attacks for a set of stochastically generated ADMM problems. Discussions on limitation and extensions are provided, along with perspectives on challenge and opportunities for fully decentralized systems.

\bibliographystyle{ieeetr}
\bibliography{EnergyBlockchain} 

\begin{thebibliography}{10}

\bibitem{10.1007/978-0-387-75462-8_6}
J.~Slay and M.~Miller, ``{Lessons Learned from the Maroochy Water Breach},'' in
  {\em Critical Infrastructure Protection} (E.~Goetz and S.~Shenoi, eds.),
  (Boston, MA), pp.~73--82, Springer US, 2008.

\bibitem{Langner2011}
R.~Langner, ``{Stuxnet: Dissecting a cyberwarfare weapon},'' {\em IEEE Security
  and Privacy}, vol.~9, no.~3, pp.~49--51, 2011.

\bibitem{Greenberg2017a}
A.~Greenberg, ``{Unprecendented Malware Targets Industrial Safety Systems in
  the Middle East},'' dec 2017.

\bibitem{Greenberg2017}
A.~Greenberg, ``{How an Entire Nation Became Russia's Test Lab for Cyberwar},''
  jun 2017.

\bibitem{Conway2016}
R.~M. Lee, M.~J. Assante, and T.~Conway, ``{Analysis of the cyber attack on the
  Ukrainian power grid},'' {\em SANS Industrial Control Systems}, p.~23, 2016.

\bibitem{Perlroth2018}
N.~Perlroth and D.~Sanger, ``{Cyberattacks Put Russian Fingers on the Switch at
  Power Plants, U.S. Says},'' mar 2018.

\bibitem{Atherton2018}
K.~Atherton, ``{It's not just elections: Russia hacked the US electric grid},''
  mar 2018.

\bibitem{EO13636}
B.~Obama, ``{Executive Order 13636- Improving Critical Infrastructure
  Cybersecurity},'' 2013.

\bibitem{PPD21}
B.~Obama, ``{Presidential Policy Directive PPD-21 -- Critical Infrastructure
  Security and Resilience},'' feb 2013.

\bibitem{ISC2015}
ISC, ``{Presidential Policy Directive 21 Implementation: An Interagency
  Security Committee White Paper},'' Tech. Rep. February, 2015.

\bibitem{Rouf2012}
I.~Rouf, H.~Mustafa, M.~Xu, and W.~Xu, ``{Neighborhood watch: Security and
  privacy analysis of automatic meter reading systems},'' {\em {\ldots}
  Communications Security}, pp.~462--473, 2012.

\bibitem{Sundaram2016}
S.~Sundaram and B.~Gharesifard, ``{Distributed Optimization Under Adversarial
  Nodes},'' pp.~1--13, 2016.

\bibitem{BlanchardEPFL2017}
P.~Blanchard, E.~{Mahdi El Mhamdi}, R.~Guerraoui, and J.~Stainer, ``{Machine
  Learning with Adversaries: Byzantine Tolerant Gradient Descent},'' {\em
  Nips2017}, no.~Nips, 2017.

\bibitem{Zhang2014}
J.~Zhang, P.~Jaipuria, A.~Chakrabortty, and A.~Hussain, ``{A Distributed
  Optimization Algorithm for Attack-Resilient Wide-Area Monitoring of Power
  Systems: Theoretical and Experimental Methods},'' in {\em Decision and Game
  Theory for Security. GameSec 2014. Lecture Notes in Computer Science, vol
  8840} (R.~Poovendran and W.~Saad, eds.), pp.~350--359, Springer, 2014.

\bibitem{Liao2018}
M.~Liao and A.~Chakrabortty, ``{Optimization Algorithms for Catching Data
  Manipulators in Power System Estimation Loops},'' {\em IEEE Transactions on
  Control Systems Technology}, pp.~1--16, 2018.

\bibitem{Xie2018}
C.~Xie, O.~Koyejo, and I.~Gupta, ``{Generalized Byzantine-tolerant SGD},''
  vol.~1, 2018.

\bibitem{Yin2018}
D.~Yin, Y.~Chen, K.~Ramchandran, and P.~Bartlett, ``{Byzantine-Robust
  Distributed Learning: Towards Optimal Statistical Rates},'' 2018.

\bibitem{Su2016a}
L.~Su and N.~Vaidya, ``{Fault-Tolerant Multi-Agent Optimization : Optimal
  Distributed Algorithms},'' in {\em Proceedings of the 2016 ACM Symposium on
  Principles of Distributed Computing}, vol.~1, (Chicago, Illinois),
  pp.~425--434, 2016.

\bibitem{Chen2018}
Y.~Chen, S.~Kar, and J.~M. Moura, ``{Resilient Distributed Estimation Through
  Adversary Detection},'' {\em IEEE Transactions on Signal Processing},
  vol.~66, no.~9, pp.~2455--2469, 2018.

\bibitem{Alistarh2018}
D.~Alistarh, Z.~Allen-Zhu, and J.~Li, ``{Byzantine Stochastic Gradient
  Descent},'' pp.~1--20, 2018.

\bibitem{Nabavi2015}
S.~Nabavi and A.~Chakrabortty, ``{An Intrusion-Resilient Distributed
  Optimization Algorithm for Modal Estimation in Power Systems},'' {\em
  Conference on Decision and Control (CDC)}, no.~Cdc, 2015.

\bibitem{Liao2016}
M.~Liao and A.~Chakrabortty, ``{A Round-Robin ADMM algorithm for identifying
  data-manipulators in power system estimation},'' {\em Proceedings of the
  American Control Conference}, vol.~2016-July, pp.~3539--3544, 2016.

\bibitem{Chen2017b}
Y.~Chen, S.~Kar, and J.~M.~F. Moura, ``{Resilient Distributed Estimation:
  Sensor Attacks},'' pp.~1--8, 2017.

\bibitem{Weldehawaryat2017}
G.~K. Weldehawaryat, P.~L. Ambassa, A.~M. Marufu, S.~D. Wolthusen, and A.~V.
  Kayem, ``{Decentralised scheduling of power consumption in micro-grids:
  Optimisation and security},'' {\em Lecture Notes in Computer Science
  (including subseries Lecture Notes in Artificial Intelligence and Lecture
  Notes in Bioinformatics)}, vol.~10166 LNCS, pp.~69--86, 2017.

\bibitem{Sundaram2016a}
S.~Sundaram and B.~Gharesifard, ``{Consensus-based distributed optimization
  with malicious nodes},'' {\em 2015 53rd Annual Allerton Conference on
  Communication, Control, and Computing, Allerton 2015}, pp.~244--249, 2016.

\bibitem{Liu2011}
Y.~Liu, P.~Ning, and M.~K. Reiter, ``{False data injection attacks against
  state estimation in electric power grids},'' {\em ACM Transactions on
  Information and System Security}, vol.~14, no.~1, pp.~1--33, 2011.

\bibitem{Li2012}
X.~Li, X.~Liang, R.~Lu, X.~Shen, X.~Lin, and H.~Zhu, ``{Securing smart grid:
  Cyber attacks, countermeasures, and challenges},'' {\em IEEE Communications
  Magazine}, vol.~50, no.~8, pp.~38--45, 2012.

\bibitem{Ozay2013}
M.~Ozay, I.~Esnaola, F.~T.~Y. Vural, S.~R. Kulkarni, and H.~V. Poor, ``{Sparse
  attack construction and state estimation in the smart grid: Centralized and
  distributed models},'' {\em IEEE Journal on Selected Areas in
  Communications}, vol.~31, no.~7, pp.~1306--1318, 2013.

\bibitem{Liang2017}
G.~Liang, J.~Zhao, F.~Luo, S.~R. Weller, and Z.~Y. Dong, ``{A Review of False
  Data Injection Attacks Against Modern Power Systems},'' {\em IEEE
  Transactions on Smart Grid}, vol.~8, no.~4, pp.~1630--1638, 2017.

\bibitem{JanLiu2018}
H.~{Jan Liu}, M.~Backes, R.~Macwan, and A.~Valdes, ``{Coordination of DERs in
  Microgrids with Cybersecure Resilient Decentralized Secondary Frequency
  Control},'' vol.~9, pp.~2670--2679, 2018.

\bibitem{Kekatos2013}
V.~Kekatos and G.~B. Giannakis, ``{Distributed robust power system state
  estimation},'' {\em IEEE Transactions on Power Systems}, vol.~28, no.~2,
  pp.~1617--1626, 2013.

\bibitem{Vukovic2014}
O.~Vukovic and G.~Dan, ``{Security of fully distributed power system state
  estimation: Detection and mitigation of data integrity attacks},'' {\em IEEE
  Journal on Selected Areas in Communications}, vol.~32, no.~7, pp.~1500--1508,
  2014.

\bibitem{Boyd2010}
S.~Boyd, ``{Distributed Optimization and Statistical Learning via the
  Alternating Direction Method of Multipliers},'' {\em Foundations and
  Trends{\textregistered} in Machine Learning}, vol.~3, no.~1, pp.~1--122,
  2010.

\bibitem{LeFloch2016a}
C.~{Le Floch}, F.~Belletti, S.~Saxena, A.~M. Bayen, and S.~Moura,
  ``{Distributed optimal charging of electric vehicles for demand response and
  load shaping},'' {\em Proceedings of the IEEE Conference on Decision and
  Control}, vol.~2016-Febru, no.~Cdc, pp.~6570--6576, 2016.

\bibitem{Dobbe2016}
R.~Dobbe, D.~Arnold, S.~Liu, D.~Callaway, and C.~Tomlin, ``{Real-Time
  Distribution Grid State Estimation with Limited Sensors and Load
  Forecasting},'' 2016.

\bibitem{Mota2012}
J.~F. Mota, J.~M. Xavier, P.~M. Aguiar, and M.~P{\"{u}}schel, ``{D-ADMM: A
  distributed algorithm for compressed sensing and other separable optimization
  problems},'' {\em ICASSP, IEEE International Conference on Acoustics, Speech
  and Signal Processing - Proceedings}, no.~2, pp.~2869--2872, 2012.

\bibitem{Sulc2014}
P.~Sulc, S.~Backhaus, and M.~Chertkov, ``{Optimal Distributed Control of
  Reactive Power Via the Alternating Direction Method of Multipliers},'' {\em
  IEEE Transactions on Energy Conversion}, vol.~29, no.~4, pp.~968--977, 2014.

\bibitem{Tsai2015}
S.~C. Tsai, Y.~H. Tseng, and T.~H. Chang, ``{Communication-efficient
  distributed demand response: A randomized ADMM approach},'' {\em IEEE
  Transactions on Smart Grid}, vol.~PP, no.~99, pp.~1--11, 2015.

\bibitem{Wang2016}
Y.~Wang, L.~Wu, and S.~Wang, ``{A Fully-Decentralized Consensus-Based ADMM
  Approach for DC-OPF With Demand Response},'' {\em IEEE Transactions on Smart
  Grid}, pp.~1--11, 2016.

\end{thebibliography}

\addtolength{\textheight}{-12.7cm}   

\appendix
\section{Analytic Solution Notes}\label{sec:appendix_analyticSolution}

The x-update step is an unconstrained quadratic program, and can be solved analytically. For clarity, we let $\gamma = Bz^k + u^k$ and proceed as:
\begin{align*}
x^{k+1} &= \text{argmin}_x \quad x^T P x + c^T x + \frac{\rho}{2}||Ax+\gamma||^2_2 \\
x^{k+1} &= \text{argmin}_x \quad x^T P x + c^T x + \frac{\rho}{2}(Ax+\gamma)^T (Ax+\gamma)\\
x^{k+1} &= \text{argmin}_x \quad x^T P x + c^T x + \frac{\rho}{2}(x^T A^T Ax + 2 \gamma^T Ax + \gamma^T \gamma)\\
0 &= \frac{\partial}{\partial x} \left( x^T P x + c^T x + \frac{\rho}{2}(x^T A^T Ax + 2 \gamma^T Ax + \gamma^T \gamma) \right)\\
0 &= 2 P x + c + \rho A^T A x + \rho \gamma^T A\\
0 &= (2P + \rho A^T A)x + c + \rho \gamma^T A \\
x &= (2P + \rho A^T A)^{-1}(-c - \rho \gamma^T A)
\end{align*}

Similarly, the z-update step can be solved analytically by letting $\mu=Ax^{k+1} - c + u^k$ and following a similar process to find:
$$
z^{k+1} = (2Q + \rho B^T B)^{-1} (-d - \rho \mu^T B)
$$

These analytic solutions are used in our implementation to avoid inaccuracies induced from a numeric solution.

\subsection{Comparison with Central Solution}
Because the problem is an unconstrained QP and entries with consensus between a subset of variables, a centralized solution can be computed by composing the cost matrices into a single quadratic problem which can be solved analytically.  This is shown here for the case where $c=0$ and $A$ and $B$ are composed as described above, but also can be computed for other $A,B$.

We break P and Q into sub-matrices dependent on the number of consensus constraints $p$, where $P_{11}, Q_{00} \in \mathbb{R}^{p \times p}$ and the other dimensions follow accordingly.

\begin{align*}
P &= \left[ \begin{array}{cc} P_{00} & P_{01}\\P_{10} & P_{11} \\ \end{array} \right]\\
Q &= \left[ \begin{array}{cc} Q_{00} & Q_{01}\\Q_{10} & Q_{11}\end{array} \right]\\
\Pi &= \left[ \begin{array}{ccc}  P_{00} & P_{01} & 0 \\ P_{10} & P_{11}+Q_{00} & Q_{10} \\
0 & Q_{10} & Q_{11} \end{array} \right]
\end{align*}

Similarly, the $c$ and $d$ vectors can be combined as 
$$
\kappa = \left[ \begin{array}{c} c_0 \\ c_1 + d_0 \\ d_1 \end{array} \right]
$$

The problem can then be expressed as an unconstrained minimization problem:
$$
\text{min}_w \quad w^T \Pi w + \kappa^T w
$$
which is solved by $w^* = -\frac{1}{2}(\Pi^T) ^{-1} \kappa$

This central solution was used only for verification.

\end{document}